%% file: writeup_lrNS.tex
\pdfoutput=1
\documentclass[final, nohypdvips]{siamart0516}
\usepackage{braket,amsfonts}

\usepackage{array}

\usepackage[caption=false]{subfig}

\usepackage{pgfplots}
\usepackage{pstool} 
\usepackage{psfrag}
\usepackage{graphicx}
\usepackage{pst-pdf}
\usepackage{comment}
\usepackage{multirow}
\usepackage{mathtools}
\usepackage{array}
\newsiamthm{claim}{Claim}
\newsiamremark{rem}{Remark}
\newsiamremark{expl}{Example}
\newsiamremark{hypothesis}{Hypothesis}
\crefname{hypothesis}{Hypothesis}{Hypotheses}

\usepackage{algorithmic}

\usepackage{graphicx,epstopdf}

\usepackage{dsfont}
\usepackage{amsmath}
\usepackage{amssymb}

\usepackage{algorithmic} % <- Preamble 
\Crefname{ALC@unique}{Line}{Lines} % <- Preamble
\Crefname{ALC@unique}{Line}{Lines}

\usepackage{amsopn}

\usepackage{xspace}
\usepackage{bold-extra}
\usepackage[most]{tcolorbox}

\colorlet{texcscolor}{blue!50!black}
\colorlet{texemcolor}{red!70!black}
\colorlet{texpreamble}{red!70!black}
\colorlet{codebackground}{black!25!white!25}

\input{definitions/commands}

\lstdefinestyle{siamlatex}{%
  style=tcblatex,
  texcsstyle=*\color{texcscolor},
  texcsstyle=[2]\color{texemcolor},
  keywordstyle=[2]\color{texemcolor},
  moretexcs={cref,Cref,maketitle,mathcal,text,headers,email,url},
}

\tcbset{%
  colframe=black!75!white!75,
  coltitle=white,
  colback=codebackground, % bottom/left side
  colbacklower=white, % top/right side
  fonttitle=\bfseries,
  arc=0pt,outer arc=0pt,
  top=1pt,bottom=1pt,left=1mm,right=1mm,middle=1mm,boxsep=1mm,
  leftrule=0.3mm,rightrule=0.3mm,toprule=0.3mm,bottomrule=0.3mm,
  listing options={style=siamlatex}
}

\newtcblisting[use counter=example]{example}[2][]{%
  title={Example~\thetcbcounter: #2},#1}

\newtcbinputlisting[use counter=example]{\examplefile}[3][]{%
  title={Example~\thetcbcounter: #2},listing file={#3},#1}

\DeclareTotalTCBox{\code}{ v O{} }
{ %fontupper=\ttfamily\color{texemcolor},
  Iupper=\ttfamily\color{black},
  nobeforeafter,
  tcbox raise base,
  colback=codebackground,colframe=white,
  top=0pt,bottom=0pt,left=0mm,right=0mm,
  leftrule=0pt,rightrule=0pt,toprule=0mm,bottomrule=0mm,
  boxsep=0.5mm,
  #2}{#1}

\patchcmd\newpage{\vfil}{}{}{}
\begin{tcbverbatimwrite}{tmp_\jobname_header.tex}
\title{A Low-rank solver for the Navier--Stokes equations with uncertain viscosity%
  \thanks{This work was supported by the U.S. Department of Energy Office of Advanced Scientific Computing Research, Applied Mathematics program under award DEC-SC0009301 and by the U.S. National Science Foundation under grant DMS1418754 and DMS1521563.}}

\author{Kookjin Lee%
  \thanks{Department of Computer Science, University of Maryland, College Park, MD 20742 (\email{klee@cs.umd.edu}).}%
  \and
  Howard C. Elman%
  \thanks{Department of Computer Science and Institute for Advanced Computer Studies, University of Maryland, College Park, MD 20742
    (\email{elman@cs.umd.edu}).}
      \and
   Bed\v{r}ich Soused\'{i}k%
  \thanks{Department of Mathematics and Statistics, University of Maryland, Baltimore County, MD, 21250 (\email{sousedik@umbc.edu}).}
}

\headers{Low-rank approximation for the Navier--Stokes equations}
{K.\ Lee, H.\ C.\ Elman, and B.\ Soused\'{i}k}
\end{tcbverbatimwrite}
\input{tmp_\jobname_header.tex}

\begin{document}
\maketitle

%% ------------------------------------------------------------------
%% ABSTRACT
%% ------------------------------------------------------------------
\begin{tcbverbatimwrite}{tmp_\jobname_abstract.tex}
\begin{abstract}
We study an iterative low-rank approximation method for the solution of the steady-state stochastic Navier--Stokes equations with uncertain viscosity. The method is based on linearization schemes using Picard and Newton iterations and stochastic finite element discretizations of the linearized problems. For computing the low-rank approximate solution, we adapt the nonlinear iterations to an inexact and low-rank variant, where the solution of the linear system at each nonlinear step is approximated by a quantity of low rank. This is achieved by using a tensor variant of the GMRES method as a solver for the linear systems. We explore the inexact low-rank nonlinear iteration with a set of benchmark problems, using a model of flow over an obstacle, under various configurations characterizing the statistical features of the uncertain viscosity, and we demonstrate its effectiveness by extensive numerical experiments.
\end{abstract}


% 35R60, 60H15: SPDE
% 65F10: Iterative methods for linear systems
\end{tcbverbatimwrite}
\input{tmp_\jobname_abstract.tex}
%% ------------------------------------------------------------------
%% END HEADER
%% ------------------------------------------------------------------

\input{introduction}

\input{stochasticGalerkin}

\input{lowrankmethods}

\input{adaptive_tolerances}
\input{numerical_results2}
\input{conclusion}

\bibliography{references}
\bibliographystyle{siamplain}

\end{document}

%% file: definitions/commands.tex
\newcommand{\truncnl}[1]{\epsilon_{\text{nl}}}
\newcommand{\truncsol}[1]{\epsilon_{\text{trunc,sol}}}
\newcommand{\truncgmres}[1]{\epsilon_{\text{gmres}}}
\newcommand{\trunccorr}[1]{\epsilon_{\text{trunc,corr}}}
\newcommand{\resnl}[1]{r_{\text{nl}}}
\newcommand{\resiter}[1]{r_{\text{gm}}}
\newcommand{\subgmres}[1]{#1_{\text{gm}}}
\newcommand{\ranksym}[1]{\alpha_{#1}}

\newcommand{\defeq}{\coloneqq}

\newcommand{\nlres}[1]{\mathcal R({#1})}
\newcommand{\vecrep}[1]{\bar #1}

\newcommand{\REV}[1]{\textcolor{blue}{#1}}

%% file: introduction.tex
\section{Introduction} \label{sec:intro} 
We are interested in the efficient computation of solutions of the steady-state Navier--Stokes equations with uncertain viscosity. Such uncertainty may arise from measurement error or uncertain ratios of multiple phases in porous media. The uncertain viscosity can be modeled as a positive random field parameterized by a set of random variables \cite{powell2012preconditioning, Sousedik2016435, tamellini2014model} and, consequently, the solution of the stochastic Navier--Stokes equations also can be modeled as a random vector field depending on the parameters associated with the viscosity (i.e., a function of the same set of random variables).  As a solution method, we consider the stochastic Galerkin method \cite{babuska2004galerkin,  ghanem2003stochastic} combined with the generalized polynomial chaos (gPC) expansion \cite{xiu2002wiener}, which provides a spectral approximation of the solution function. The stochastic Galerkin method results in a coupled algebraic system of equations.  There has been considerable progress in development of solvers for these systems \cite{elman2007solving, powell2009block,rosseel2010iterative, sousedik2014hierarchical,ullmann2012efficient}, although costs may be high when the global system becomes large.

One way to address this issue is thorough use of tensor \textit{Krylov subspace} methods, which operate in tensor format and reduce the costs of matrix operations by exploiting a Kronecker-product structure of  system matrices. Variants of this approach have been developed for the Richardson iteration \cite{kressner2011low, matthies2012solving}, the conjugate gradient method \cite{kressner2011low}, the BiCGstab method \cite{kressner2011low}, the minimum residual method \cite{stoll2015low}, and the general minimum residual (GMRES) method \cite{ballani2013projection}. Efficiencies are also obtained from the fact that solutions can often be well approximated by low-rank objects. These ideas have been shown to reduce costs for solving steady \cite{lee2016preconditioned, matthies2012solving} and unsteady stochastic diffusion equations \cite{benner2015low}.

In this study, we adapt the low-rank approximation scheme to a solver for the systems of nonlinear equations obtained from the stochastic Galerkin discretization of the stochastic Navier--Stokes equations. In particular, we consider a low-rank variant of linearization schemes based on Picard and Newton iteration, where the solution of the nonlinear system is computed by solving a sequence of linearized systems using a low-rank variant of the GMRES method (lrGMRES) \cite{ballani2013projection} in combination with inexact nonlinear iteration \cite{inexactNewton1982}. 

We base our development of the stochastic Galerkin formulation of the stochastic Navier--Stokes equations  on ideas from \cite{powell2012preconditioning,Sousedik2016435}. In particular, we consider a random viscosity affinely dependent on a set of random variables as suggested in \cite{powell2012preconditioning} (and in \cite{Sousedik2016435}, which considers a gPC approximation of the lognormally distributed viscosity). The stochastic Galerkin formulation of the stochastic Navier--Stokes equations is also considered in \cite{benner2017lrNS}, which studies an optimal control problem constrained by the stochastic Navier--Stokes problem and computes an approximate solution using a low-rank tensor-train decomposition \cite{oseledets2011tensor}. Related work \cite{tamellini2014model} extends a Proper Generalized Decomposition method \cite{nouy2007generalized} for the stochastic Navier--Stokes equations, where a low-rank approximate solution is built from successively computing rank-one approximations. See the book \cite{le2010spectral} for an overview and other spectral approximation approaches for models of computational fluid dynamics.

An outline of the paper is as follows. In section \ref{sec:problem}, we review the stochastic Navier--Stokes equations and their discrete Galerkin formulations. In section \ref{sec:lr_method}, we present an iterative low-rank approximation method for solutions of the discretized stochastic Navier--Stokes problems. In section \ref{sec:adap_tol}, we introduce an efficient variant of the inexact Newton method, which solves linear systems arising in nonlinear iteration using low-rank format. We follow a hybrid approach, which employs several steps of Picard iteration followed by Newton iteration. In section \ref{sec:results}, we examine the performance of the proposed method on a set of benchmark problems that model the flow over an obstacle. Finally, in section \ref{sec:conclusion}, we draw some conclusions.

%% file: stochasticGalerkin.tex
\section{Stochastic Navier--Stokes equations}\label{sec:problem}
Consider the stochastic Navier--Stokes equations: Find velocity $\vec u(x, \xi)$ and pressure $p(x, \xi)$ such that
\begin{align}\label{eq:NS_strong}
\begin{split}
-\nu(x, \xi) \nabla^2 \vec u(x, \xi) + (\vec u(x, \xi) \cdot \nabla)\vec u(x, \xi) + \nabla p(x, \xi)  &= \vec f(x, \xi),\\
\nabla \cdot \vec u(x, \xi) &= 0,
\end{split}
\end{align}
in $D \times \Gamma$, with a boundary conditions
\begin{alignat*}{2}
\vec u (x,\xi) &= \vec g(x,\xi), \quad && \text{ on } \partial D_{\text{Dir}},\\
\nu(x,\xi) \nabla \vec u(x,\xi) \cdot \vec n - p(x,\xi) \vec n (x,\xi)  &= \vec 0, \quad &&\text{ on } \partial D_{\text{Neu}},
\end{alignat*}
where $\partial D = \partial D_{\text{Dir}} \, \cup \, \partial D_{\text{Neu}}$. %Here, $\nu(x,\xi)$ is a random viscosity, 
The stochasticity of the equation \eqref{eq:NS_strong} stems from the random viscosity $\nu(x, \xi)$, which is modeled as a positive random field parameterized by a set of independent, identically distributed random variables $\xi = \{\xi_1,\ldots,\xi_{n_\nu}\}$. The random variables comprising $\xi$ are defined on a probability space $(\Omega, \mathcal F, P)$ such that $\xi: \Omega \rightarrow \Gamma \subset \mathbb R^{n_\nu}$, where $\Omega$ is a sample space, $\mathcal F$ is a $\sigma$-algebra on $\Omega$, and $P$ is a probability measure on $\Omega$. The joint probability density function of $\xi$ is denoted by $\rho(\xi)$ and the expected value of a random function $v(\xi)$ on $\Gamma$ is then $\langle v \rangle_{\rho} = \mathbb E[v ] \equiv \int_\Gamma v(\xi) \rho(\xi) d\xi$.

For the random viscosity, we consider a random field that has affine dependence on the random variables $\xi$,
\begin{align}\label{eq:rf_visc}
\nu(x,\xi) \equiv \nu_0 + \sigma_\nu \sum_{k=1}^{n_\nu} \nu_k(x) \xi_k,
\end{align}
where $\{\nu_0,\sigma_\nu^2 \}$ are the mean and the variance of the random field $\nu(x,\xi)$. We will also refer to the coefficient of variation ($CoV$), the relative size of the standard deviation with respect to the mean,
\begin{align}\label{eq:cov_def}
CoV \equiv \frac{\sigma_\nu}{\nu_0}.
\end{align}
The random viscosity leads to the random Reynolds number
\begin{align}\label{eq:reynold_def}
\text{Re}(\xi) \equiv \frac{UL}{\nu(\xi)},
\end{align}
where $U$ is the characteristic velocity and $L$ is the characteristic length. We denote the Reynolds number associated with the mean viscostiy by $\text{Re}_0 = \frac{UL}{ \nu_0}$. In this study, we ensure that the viscosity \eqref{eq:rf_visc} has positive values by controlling $CoV$ and only consider small enough  $\text{Re}_0$ so that the flow problem has a unique solution. 

\subsection{Stochastic Galerkin method}
In the stochastic Galerkin method, a mixed variational formulation of \eqref{eq:NS_strong} can be obtained by employing Galerkin orthogonality: Find $(\vec u, p) \in (V_E, Q_D) \otimes L^2(\Gamma)$ such that 
\begin{alignat}{2} 
\left \langle \int_D \nu \nabla \vec u : \nabla \vec v + (\vec u \cdot \nabla \vec u) \vec v -  p(\nabla \cdot \vec v) \right \rangle_\rho &= \left \langle \int_D \vec f \cdot  \vec v  \right \rangle_\rho, \quad &&\forall \vec v \in V_D \otimes L^2(\Gamma),\label{eq:NS_weak1}\\
\left\langle \int_D q(\nabla \cdot \vec u ) \right \rangle_\rho &= 0, &&\forall q \in Q_D \otimes L^2(\Gamma).\label{eq:NS_weak2}
\end{alignat} 
The velocity solution and test spaces are $V_E = \{ \vec u \in \mathcal H^1(D)^2 | \vec u = \vec g \text{ on } \partial D_{\text{Dir}}\}$ and $V_D = \{ \vec v \in \mathcal H^1(D)^2 | \vec v = \vec 0 \text{ on } \partial D_{\text{Dir}}\}$, where $\mathcal H^1(D)$ refers to the Sobolev space of functions with derivatives in $L^2(D)$, for the pressure solution, $Q_D =  L^2(D)$, and $L^2(\Gamma)$ is a Hilbert space equipped with an inner product
\begin{equation*}\label{eq:inner_prod}
\langle u, v \rangle_\rho \equiv \int_\Gamma u(\xi) v(\xi) \rho(\xi) d\xi.
\end{equation*} 
 %The solution space $(V_D , Q_D )$ satisfies an inf-sup condition.
The solution of the variational formulation \eqref{eq:NS_weak1}--\eqref{eq:NS_weak2} satisfies
\begin{equation}\label{eq:nonlinear_system}
\nlres{\vec u, p; \vec v, q} = 0, \qquad \forall \vec v \in V_D\otimes L^2(\Gamma),\, \forall q\in Q_D\otimes L^2(\Gamma),
\end{equation}
where $\nlres{\vec u, p; \vec v,q}$ is a nonlinear residual
\begin{align}\label{eq:residual_def}
\nlres{\vec u, p; \vec v,q} \equiv \begin{bmatrix} \langle  \int_D \vec f \cdot  \vec v - \nu \nabla \vec u : \nabla \vec v +  (\vec u \cdot \nabla \vec u) \vec v - \int_D p(\nabla \cdot \vec v)  \rangle_\rho\\ 
\langle - \int_D q (\nabla \cdot \vec u) \rangle_\rho
 \end{bmatrix}.
\end{align}
To compute the solution of the nonlinear equation \eqref{eq:nonlinear_system}, we employ linearization techniques based on either Picard iteration or Newton iteration \cite{elman2014finite}. 
Replacing $(\vec u, p)$ of \eqref{eq:NS_weak1}--\eqref{eq:NS_weak2} with $(\vec u + \delta \vec u, p + \delta p)$  and neglecting the quadratic term $c(\delta \vec u; \delta \vec u, \vec v)$, where $c( \vec z; \vec u, \vec v) \equiv \int_D (\vec z \cdot \nabla \vec u) \cdot \vec v$, gives
\begin{align}\label{eq:jacobian}
\begin{bmatrix} \langle \int_D \nu \nabla \delta \vec u : \nabla \vec v  + c(\delta \vec u; \vec u, \vec v) + c(\vec u; \delta \vec u,\vec v) - \int_D \delta p (\nabla \cdot \vec v) \rangle_\rho \\ 
\langle \int_D q(\nabla \cdot \delta \vec u) \rangle_\rho
 \end{bmatrix} = \mathcal R (\vec u, p; \vec v, q).
\end{align}
In Newton iteration, the $(n+1)$st iterate $(\vec u^{n+1}, p^{n+1})$ is computed by taking $\vec u = \vec u^n$, $p=p^n$ in \eqref{eq:jacobian}, solving \eqref{eq:jacobian} for $(\delta \vec u^{n}, \delta p^{n})$,  and updating 
\begin{align*}
\vec u^{n+1} \defeq \vec u^n + \delta \vec u^n,\quad p^{n+1} \defeq p^n + \delta p^n.
\end{align*}
In Picard iteration, the term $c(\delta \vec u; \vec u, \vec v)$ is omitted from the linearized form \eqref{eq:jacobian}.

\subsection{Discrete stochastic Galerkin system}\label{sec:sg_disc}
To obtain a discrete system, the velocity $\vec u(x,\xi)$ and the pressure $p(x,\xi)$ are approximated by a generalized polynomial chaos expansion \cite{xiu2002wiener}:
\begin{equation}
\vec u(x,\xi) \equiv \sum_{i=1}^{n_\xi}  \vec u_i(x) \psi_i(\xi), \quad p(x,\xi) \equiv \sum_{i=1}^{n_\xi}  p_i(x) \psi_i(\xi),
\label{eq:gpc_sol}
\end{equation}
where $\{\psi_i(\xi)\}_{i=1}^{n_\xi}$ is a set of $n_\nu$-variate orthogonal polynomials (i.e., $\langle \psi_i \psi_j \rangle_\rho=0$ if $i\neq j$). This set of orthogonal polynomials gives rise to a finite-dimensional approximation space $S = \text{span}(\{\psi_i(\xi)\}_{i=1}^{n_\xi}) \subset L^2(\Gamma)$. For spatial discretization, a div-stable mixed finite element method \cite{elman2014finite} is considered, the Taylor-Hood element consisting of biquadratic velocities and bilinear pressure. Basis sets for the velocity space $V_E^h$ and the pressure space $Q_D^h$ are denoted by $\left\{ \begin{bmatrix} \phi_i(x) \\ 0 \end{bmatrix}, \begin{bmatrix} 0\\ \phi_i(x) \end{bmatrix} \right\}_{i=1}^{n_u}$ and $\{\varphi_i(x)\}_{i=1}^{n_p}$, respectively. Then the fully discrete version of \eqref{eq:gpc_sol} can be written as 
\begin{equation} \label{eq:sg_sol}
\vec { u}(x,\xi) = \begin{bmatrix} \vec u^x(x,\xi) \\ \vec u^y(x,\xi) \end{bmatrix} \equiv \begin{bmatrix}  \sum_{i=1}^{n_\xi} \sum_{j=1}^{n_u}   u^x_{ij} \phi_j(x) \psi_i(\xi) \\   \sum_{i=1}^{n_\xi} \sum_{j=1}^{n_u}  u^y_{ij} \phi_j(x) \psi_i(\xi)  \end{bmatrix}, \quad  p(x,\xi) \equiv  \sum_{i=1}^{n_\xi} \sum_{j=1}^{n_p}   p_{ij}  \varphi_j(x) \psi_i(\xi).
\end{equation}
Let us introduce a vector notation for the coefficients, $\vecrep{u}_i^x \equiv [u_{i1}^x, \ldots, u_{in_u}^x]^T \in \mathbb{R}^{n_u}$, $\vecrep{u}_i^y \equiv [u_{i1}^y, \ldots, u_{in_u}^y]^T \in \mathbb{R}^{n_u}$, and $\vecrep{p}_i \equiv  [p_{i1}, \ldots, p_{in_p}]^T \in \mathbb{R}^{n_p}$ for $i=1,\ldots,n_\xi$, which, for each gPC index i, groups the horizontal velocity coefficients together followed by the vertical velocity coefficients, and then by the pressure coefficients, giving a vector
\begin{equation}\label{eq:gPC_coeff_ith}
\bar u_i = [(\bar u^x_i)^T,  (\bar u^y_i)^T, p_i^T]^T.
\end{equation}
Taking $\nu(x,\xi)$ from \eqref{eq:rf_visc} and replacing $\vec u(x,\xi)$, $p(x,\xi)$ in \eqref{eq:jacobian} with their discrete approximations \eqref{eq:sg_sol} yields a system of linear equations of order $(2n_u + n_p) n_\xi$.
The coefficient matrix has a Kronecker-product structure, 
\begin{align} \label{eq:jacobian_matrix}
%J^n =  G_0 \otimes \mathcal F^n_0 + \sum_{l=1}^{n_\xi} G_l \otimes \mathcal F^n_l,
J \equiv  G_1 \otimes \mathcal F_1 + \sum_{l=2}^{n_\xi} G_l \otimes \mathcal F_l,
\end{align}
where $G_l$ refers to the $l$th ``stochastic matrix'' 
\begin{align*} \label{eq:stochastic_matrix}
[G_l]_{ij} = \langle \psi_l \psi_i \psi_j \rangle_\rho, \quad l = 1,\ldots, n_\xi
\end{align*}
with $\psi_1(\xi) = 1$, $\psi_i(\xi) = \xi_{i-1}$ for $i=2,\ldots,n_\nu+1$ and
\begin{align*} \label{eq:newton_deterministic}
\mathcal F_1 \equiv \begin{bmatrix} F_1 & B^T \\ B & 0 \end{bmatrix}, \quad \mathcal F_l \equiv \begin{bmatrix} F_l & 0 \\ 0 & 0 \end{bmatrix}, \quad l=2,\ldots,n_\xi
\end{align*}
with
$F_l \equiv A_l + N_l + W_l$ for the Newton iteration and $F_l \equiv A_l + N_l$  for the Picard iteration. We refer to the matrix of \eqref{eq:jacobian_matrix} derived from the Newton iteration as the \textit{Jacobian} matrix, and that derived from the Picard iteration as the \textit{Oseen} matrix, denoted by $J_N$ and $J_P$, respectively. Here, $A_l$ is the $l$th symmetric matrix defined as
\begin{align}\label{eq:laplace_matrix}
[A_l]_{ij} \equiv \int_D \nu_{l-1} (x) (\nabla \phi_i : \nabla \phi_j ), \quad l=1,\ldots,n_\nu+1,
\end{align}
$N_l = N(\vec u_l(x))$ and $W_l = W(\vec u_l(x))$ are, respectively, the $l$th vector-convection matrix and the $l$th Newton derivative matrix with $\vec u^n_l(x)$ from the $l$th term of \eqref{eq:gpc_sol},
\begin{align*}
[N_l]_{ij} = [N(\vec u_l(x))]_{ij} &\equiv \int_D (\vec u_l(x) \cdot \nabla \phi_j(x)) \cdot \phi_i(x), \quad l=1,\ldots,n_\xi,\\
[W_l]_{ij} = [W(\vec u_l(x))]_{ij} &\equiv \int_D ( \phi_j(x) \cdot \nabla \vec u_l(x)) \cdot \phi_i(x), \quad l=1,\ldots,n_\xi,
\end{align*}
and $B$ is the divergence matrix,
\begin{align} \label{eq:div_matrix}
[B]_{ij} \equiv \int_D \varphi_j (\nabla \cdot \phi_i).
\end{align}
If the number of gPC polynomial terms in \eqref{eq:sg_sol} is larger than the number of terms in \eqref{eq:rf_visc} (i.e., $n_\xi > n_\nu+1$), we simply set $\{A_l\}_{l={n_\nu}+2}^{n_\xi}$ as matrices containing only zeros so that $\mathcal F_l = N_l + W_l$ for $l=n_\nu+2,\ldots,n_\xi$.
 
A discrete version of \eqref{eq:residual_def} can be derived in a similar way,
\begin{align}\label{eq:rhs_vector}
\vecrep{r} \defeq \vecrep{y} - \left(G_1 \otimes \mathcal P_1 +  \sum_{l=2}^{n_\xi} G_l \otimes \mathcal P_l \right) \vecrep{u}
\end{align}
where $\vecrep{u} \defeq [\bar u_1^T \ldots \bar u_{n_\xi}^T]^T \in \mathbb{R}^{(2n_u + n_p) n_\xi}$ with $\bar u_i$ as in \eqref{eq:gPC_coeff_ith}, $\vecrep{y}$ is the right-hand side determined from the forcing function and Dirichlet boundary data, and 
\begin{align*}
\mathcal P_1 \equiv \begin{bmatrix} A_1 + N_1 & B^T \\ B & 0 \end{bmatrix}, \quad \mathcal P_l \equiv \begin{bmatrix} A_l + N_l & 0 \\ 0 & 0 \end{bmatrix} \quad l=2,\ldots,n_\xi. 
\end{align*}
The system of linear equations arising at the $n$th nonlinear iteration is
\begin{align}\label{eq:jacobian_system}
J^n \delta \vecrep{u}^n =- \vecrep{r}^n,
\end{align}
where the matrix $J^n$ from  \eqref{eq:jacobian_matrix} and the residual $\vecrep{r}^n$ from \eqref{eq:rhs_vector} each evaluated at the $n$th iterate $\vecrep{u}^n$, and the update $\delta \bar u^n$ is computed by solving \eqref{eq:jacobian_system}. The order of the system $(2n_u + n_p)n_\xi$ grows fast as the number of random variables used to parameterize the random viscosity increases. Even for a moderate-dimensional stochastic Navier--Stokes problem, solving a sequence of linear systems of order $(2n_u + n_p)n_\xi$ can be computationally prohibitive. To address this issue, we present an efficient variant of Newton--Krylov methods in the following sections.

%% file: lowrankmethods.tex
\section{Low-rank Newton--Krylov method}\label{sec:lr_method}
In this section, we outline the formalism in which the solutions to \eqref{eq:rhs_vector} and \eqref{eq:jacobian_system} can be efficiently approximated by low-rank objects while not losing much accuracy and we show how solvers are adjusted within this formalism.

Before presenting these ideas, we describe the nonlinear iteration. We consider a hybrid strategy. An initial approximation for the nonlinear solution is computed by solving the parameterized Stokes equations, 
\begin{equation*}
\begin{split}
-\nu(x, \xi) \nabla^2 \vec u(x, \xi) + \nabla p(x, \xi)  &= \vec f(x, \xi),\\
\nabla \cdot \vec u(x, \xi) &= 0\REV{.}
\end{split}
\end{equation*}
The discrete Stokes operator, which is obtained from the stochastic Galerkin discretization as shown in section \ref{sec:sg_disc}, is
\begin{align}\label{eq:stokes_matrix}
\left(G_1 \otimes \mathcal S_1 + \sum_{l=2}^{n_\nu+1} G_l \otimes \mathcal S_l \right) \bar u_{\text{st}} = b_{\text{st}},
\end{align} 
where 
\begin{align*}
\mathcal S_1 = \begin{bmatrix} A_1 & B^T \\ B & 0 \end{bmatrix}, \qquad \mathcal S_l = \begin{bmatrix} A_l & 0 \\ 0 & 0 \end{bmatrix}, \quad l = 2,\ldots,n_\nu+1,
\end{align*}
with $\{A_l\}_{l=1}^{n_\nu+1}$ defined in \eqref{eq:laplace_matrix} and $B$ defined in \eqref{eq:div_matrix}.
After this initial computation, updates to the solution are computed by first solving $m_p$ Picard systems with coefficient matrix $J_P$ and then using Newton's method with coefficient matrix $J_N$ to compute the solution.

\begin{algorithm}[H]
\caption{Solution methods}
\label{alg:abs_proc}
\begin{algorithmic}[1]
\STATE compute an approximate solution of $A_{\text{st}} \bar u_{\text{st}} = b_{\text{st}}$ in \eqref{eq:stokes_matrix}% using Algorithm \ref{alg:lrp_alg}
\STATE{set an initial guess for the Navier--Stokes problem $\bar u^0 \defeq \bar u_{\text{st}}$}\\
\FOR[Picard iteration]{ $k = 0,\ldots,m_p-1$}
\STATE solve $J_{\text{P}}^k\, \delta \vecrep{u}^k = -\vecrep{r}^k$ 
\STATE update $\bar u^{k+1} \defeq  \bar u^k + \delta \bar u^k$\label{ln:trunc1}
\ENDFOR\\
\WHILE[Newton iteration]{ $k < m_n$ and $\| \vecrep{r}^k  \|_2 > \truncnl{} \| \bar r ^0\|_2$ }
\STATE solve $J_{\text{N}}^k\, \delta \bar u^k = -\vecrep{r}^k$ 
\STATE update $\bar u^{k+1} \defeq  \bar u^k + \delta \bar u^k$ \label{ln:trunc2}
\ENDWHILE
\end{algorithmic}
\end{algorithm}

\subsection{Approximation in low rank}
We now develop a low-rank variant of Algorithm \ref{alg:abs_proc}. Let us begin by introducing some concepts to define the rank of computed quantities. Let $X =[\vecrep{x}_1, \cdots, \vecrep{x}_{n_2}] \in \mathbb{R}^{n_1 \times n_2}$ and $\vecrep{x}  = [\vecrep{x}_1^T, \cdots, \vecrep{x}_{n_2}^T]^T \in \mathbb{R}^{n_1n_2}$, where $\vecrep{x}_i \in \mathbb{R}^{n_1}$ for $i=1,\ldots,n_2$. That is, $\vecrep{x}$ can be constructed by rearranging the elements of $X$, and vice versa. Suppose $X$ has rank $\ranksym{x}$. Then two mathematically equivalent expressions for $X$ and $\vecrep{x}$ are given by
\begin{equation}\label{eq:tensor_rep}
X = YZ^T = \sum_{i=1}^{\ranksym{\vecrep{x}}} \vecrep{y}_i \vecrep{z}_i^T \quad  \Leftrightarrow \quad \vecrep{x} =  \sum_{i=1}^{\ranksym{\vecrep{x}}} \vecrep{z}_i \otimes \vecrep{y}_i,
\end{equation}
where $Y\equiv[\vecrep{y}_1, \cdots, \vecrep{y}_{\ranksym{\vecrep{x}}}] \in \mathbb{R}^{n_1 \times \ranksym{\vecrep{x}}}$, $Z\equiv[\vecrep{z}_1, \cdots, \vecrep{z}_{\ranksym{\vecrep{x}}}] \in \mathbb{R}^{n_2 \times \ranksym{\vecrep{x}}}$ with $\vecrep{y}_i \in \mathbb{R}^{n_1}$,  $\vecrep{z}_i \in  \mathbb{R}^{n_2}$ for $i=1,\ldots,\ranksym{\vecrep{x}}$. The representation of $X$ and its rank is standard matrix notation; we also use $\ranksym{x}$ to refer to the rank of the corresponding vector $\vecrep{x}$.

With this definition of rank, our goal is to inexpensively find a low-rank approximate solution $\vecrep{u}^k$ satisfying $\| \vecrep{r}^k  \|_2 \leq \truncnl{} \| \bar r ^0\|_2$ for small enough $\truncnl{}$. To achieve this goal, we approximate updates $\{\delta \vecrep{u}^k\}$ in low-rank using a low-rank variant of GMRES method, which exploits the Kronecker product structure in the system matrix as in \eqref{eq:jacobian_matrix} and \eqref{eq:stokes_matrix}. In the following section, we present the solutions $\vecrep{u}$ (and $\delta \vecrep{u}$) in the formats of \eqref{eq:tensor_rep} together with matrix and vector operations that are essential for developing the low-rank GMRES method.

\subsection{Solution coefficients in Kronecker-product form}\label{sec:matricized_sol}
We seek separate low-rank approximations of the horizontal and vertical velocity solutions and the pressure solution. With the representation shown in \eqref{eq:tensor_rep}, the solution coefficient vector $\bar u \in \mathbb{R}^{(2n_u + n_p) n_\xi}$, which consists of the coefficients of the velocity solution and the pressure solution \eqref{eq:sg_sol}, has an equivalent representation $U \in \mathbb{R}^{(2n_u + n_p) \times n_\xi}$. The matricized solution coefficients $U = [{U^x}^T, {U^y}^T, P^T]^T$ where $U^x = [\bar u^x_1,\ldots,\bar u^x_{n_\xi}], U^y = [\bar u^y_1,\ldots,\bar u^y_{n_\xi}] \in \mathbb{R}^{n_u \times n_\xi}$ and the pressure solution $P= [\bar p_1,\ldots,\bar p_{n_\xi}] \in \mathbb{R}^{n_p \times n_\xi}$.
The components admit the following representations:
\begin{align}
U^x &= \sum_{i=1}^{\ranksym{\bar u^x}} \vecrep{v}^x_i (\vecrep{w}^x_i)^T = V^x(W^x)^T  \quad \Leftrightarrow \quad \bar u^{x} = \sum_{i=1}^{\ranksym{\bar u^x}} \vecrep{w}^{x}_i \otimes \vecrep{v}^{x}_i, \label{eq:ux_tensor}\\
U^y &= \sum_{i=1}^{\ranksym{\bar u^y}} \vecrep{v}^y_i (\vecrep{w}^y_i)^T = V^y(W^y)^T  \quad \Leftrightarrow \quad \bar u^{y} = \sum_{i=1}^{\ranksym{\bar u^y}} \vecrep{w}^{y}_i \otimes \vecrep{v}^{y}_i,\label{eq:uy_tensor}\\
P &= \sum_{i=1}^{\ranksym{\bar p}} \vecrep{v}^p_i (\vecrep{w}^p_i)^T= V^p(W^p)^T \quad \Leftrightarrow \quad \bar p = \sum_{i=1}^{\ranksym{\bar p}} \vecrep{w}^{p}_i  \otimes \vecrep{v}^{p}_i, \label{eq:p_tensor}
\end{align}
where $V^x = [\vecrep{v}^{x}_1 \ldots \vecrep{v}^{x}_{\ranksym{\bar u^x}}]$, $W^x = [\vecrep{w}^{x}_1 \ldots \vecrep{w}^{x}_{\ranksym{\bar u^x}}]$, $\ranksym{\bar u^x}$ is the rank of $\bar u^x$ and $U^x$, and the same interpretation can be applied to $\bar u^y$ and $\bar p$.

\subsubsection{Matrix operations}\label{sec:matrix_op}
In this section, we introduce essential matrix operations used by the low-rank GMRES methods, using the representations shown in~\eqref{eq:ux_tensor}--\eqref{eq:p_tensor}. First, consider  the matrix-vector product with the Jacobian system matrix  \eqref{eq:jacobian_matrix} and vectors \eqref{eq:ux_tensor}--\eqref{eq:p_tensor}, 
\begin{align}\label{eq:mvp_tensor}
J^n \bar u^n = \left(\sum_{l=1}^{n_\xi} G_l \otimes \mathcal F_l^n\right) \bar u^n,
\end{align}
where 
\begin{align*}
\mathcal F^n_l \!=\! \begin{bmatrix} A_l^{xx} + N_l^n + W_l^{xx,n} & W_l^{xy,n} & {B^x}{}^T \\
 W_l^{yx,n} & A_l^{yy} + N_l^n + W_l^{yy,n}  & B^y{}^T \\ 
B^x & B^y & 0 \end{bmatrix} \!=\! \begin{bmatrix} \mathcal F_l^{xx, n} & \mathcal F_l^{xy,n} & B^x{}^T \\
\mathcal F_l^{yx,n} & \mathcal F_l^{yy,n}   & B^y{}^T \\ 
B^x & B^y & 0 \end{bmatrix}
\end{align*}
with $\mathcal F^{xx,n}_l, \mathcal F^{xy,n}_l, \mathcal F^{yx,n}_l, \mathcal F^{yy,n}_l \in \mathbb{R}^{n_u \times n_u}$ and $B^x, B^y \in \mathbb{R}^{n_p \times n_u}$. The expression \eqref{eq:mvp_tensor} has the equivalent matricized form $\sum_{l=1}^{n_\xi}  \mathcal F_l^n U^n G_l^T$ where the $l$th-term is evaluated as
\begin{align}\label{eq:mvp_detail}
 \mathcal F_l^n U^n G_l^T \!=\! \! \begin{bmatrix} \mathcal F_l^{xx, n} V^{x,n} (G_lW^{x,n})^T \! + \! F_l^{xy,n}  V^{y,n} (G_lW^{y,n})^T \!+\! {B^x}^T V^{p,n} (G_lW^{p,n})^T \\
\mathcal F_l^{yx, n} V^{x,n} (G_lW^{x,n})^T\! +\!  F_l^{yy,n}  V^{y,n} (G_lW^{y,n})^T \!+\! {B^y}^T V^{p,n} (G_lW^{p,n})^T \\
B^x V^{x,n} ( G_lW^{x,n})^T + B^y  V^{y,n} ( G_lW^{y,n})^T
 \end{bmatrix}.
\end{align}
Equivalently, in the Kronecker-product structure, the matrix-vector product \eqref{eq:mvp_detail} updates each set of solution coefficients as follows: 
\begin{alignat}{2}
&\sum_{l=1}^{n_\xi}  (G_l \otimes \mathcal F_l^{xx,n}) \bar u^{x,n} + (G_l \otimes \mathcal F_l^{xy,n}) \bar u^{y,n} + (G_l \otimes B^{x}{}^T) \bar p^n, \quad &&(x\text{-velocity}), \label{eq:mvp_x}\\ 
&\sum_{l=1}^{n_\xi}  (G_l \otimes \mathcal F_l^{yx,n}) \bar u^{x,n} + (G_l \otimes \mathcal F_l^{yy,n}) \bar u^{y,n}  + (G_l \otimes B^{y}{}^T) \bar p^n,  &&(y\text{-velocity})\label{eq:mvp_y} \\
&\sum_{l=1}^{n_\xi} (G_l \otimes B^{x}) \bar u^{x,n} +(G_l \otimes B^{y}) \bar u^{y,n},  &&(\text{pressure})\label{eq:mvp_p}
\end{alignat}
where each matrix-vector product can be performed by exploiting the Kronecker-product structure, for example,
\begin{align}\label{eq:mvp_split2}
\begin{split}
\sum_{l=1}^{n_\xi}   (G_l \otimes \mathcal F_l^{xx,n}) \bar u^{x,n} = \sum_{l=1}^{n_\xi}   G_l \otimes \mathcal F_l^{xx,n} \sum_{i=1}^{\ranksym{\bar u^x}} w^{x}_i \otimes v^{x}_i = \sum_{l=1}^{n_\xi}\sum_{i=1}^{\ranksym{\bar u^x}}   G_l w^{x}_i \otimes \mathcal F_l^{xx,n}v^{x}_i.
\end{split}
\end{align}
The matrix-vector product shown in \eqref{eq:mvp_x}--\eqref{eq:mvp_p} requires $O( 2n_u + n_p + n_\xi)$ flops, whereas \eqref{eq:mvp_tensor} requires $O( (2n_u + n_p)n_\xi)$ flops. Thus, as the problem size grows, the additive form of the latter count grows much less rapidly than the multiplicative form for \eqref{eq:mvp_tensor}.

The addition of two vectors $\bar u^x$ and $\bar u^y$ can also be efficiently performed in the Kronecker-product structure, 
\begin{equation}\label{eq:sum_tensor}
\bar u^x + \bar u^y = \sum_{i=1}^{\ranksym{\bar u^x}} w_i^x \otimes v_i^x + \sum_{i=1}^{\ranksym{\bar u^y}} w_i^y \otimes v_i^y = \sum_{i=1}^{\ranksym{\bar u^x} + \ranksym{\bar u^y}} \hat v_i + \hat w_i,
\end{equation}
where $\hat v_i = v_i^x,\, \hat w_i = w_i^x$ for $i=1,\ldots,\ranksym{\bar u^x}$, and $\hat v_i = v_i^y,\, \hat w_i = w_i^y$ for $i=\ranksym{\bar u^x}+1,\ldots,\ranksym{\bar u^x}+\ranksym{\bar u^y}$.

Inner products can be performed with similar efficiencies. Consider two vectors $\bar x_1$ and $\bar x_2$, whose matricized representations are 
\begin{align}
X_1 = \begin{bmatrix}  Y_{11}Z_{11}^T \\ Y_{12}Z_{12}^T \\ Y_{13}Z_{13}^T\end{bmatrix}, \qquad X_2 = \begin{bmatrix}  Y_{21}Z_{21}^T \\ Y_{22}Z_{22}^T \\ Y_{23}Z_{23}^T\end{bmatrix}.
\end{align}
Then the Euclidean inner product between $x_1$ and $x_2$ can be evaluated as
\begin{align*}\label{eq:inner_prod_tensor}
\bar x_1^T \bar x_2 = \text{trace} ( (Y_{11}Z_{11}^T)^T Y_{21}Z_{21}^T ) +  \text{trace} ( (Y_{12}Z_{12}^T)^T Y_{22}Z_{22}^T ) +  \text{trace} ( (Y_{13}Z_{13}^T)^T Y_{23}Z_{23}^T ),
\end{align*}
where trace($X$) is defined as a sum of the diagonal entries of the matrix $X$.

Although the matrix-vector product and the sum, as described in \eqref{eq:mvp_split2} and \eqref{eq:sum_tensor}, can be performed efficiently, the results of \eqref{eq:mvp_split2} and \eqref{eq:sum_tensor} are represented by $n_\xi \ranksym{\bar u^x}$ terms and $\ranksym{\bar u^x}+\ranksym{\bar u^y}$ terms, respectively, which typically causes the ranks of the computed quantities to be higher than the inputs for the computations and potentially undermines the efficiency of the solution method. To resolve this issue, a truncation operator will be used to modify the result of matrix-vector products and sums and force the ranks of quantities used to be small.

\subsubsection{Truncation of $U^{x,n}$, $U^{y,n}$ and $P^n$} \label{sec:trunc}
We now explain the details of the truncation. Consider the velocity and the pressure represented in a matrix form as in \eqref{eq:ux_tensor}--\eqref{eq:p_tensor}. The best $\ranksym{}$-rank approximation of a matrix can be found by using the singular value decomposition (SVD) \cite{kressner2011low, matthies2012solving}.  Here, we define a truncation operator for a given matrix $U = VW^T$ whose rank is $\ranksym{U}$,
\begin{align*}
\mathcal T_{\epsilon_{\text{trunc}}}&:  U \rightarrow \tilde U,
\end{align*}
where the rank of $U$ is larger than the rank of $\tilde U$ (i.e., $\ranksym{U} \gg \ranksym{\tilde U})$. The truncation operator $\mathcal T_{\epsilon_{\text{trunc}}}$ compresses $U$ to $\tilde U$ such that $\| \tilde  U - U \|_F \leq \epsilon_{\text{trunc}}\|  U\|_F$ where $\|\cdot \|_F$ is the Frobenius norm. To achieve this goal, the singular value decomposition of $U$ can be computed (i.e., $U = \hat V D \tilde W^T$ where $D=\text{diag}(d_1,\ldots,d_n)$ is the diagonal matrix of singular values). Letting $\{\hat v_i\}$ and $\{\tilde w_i\}$ denote the singular vectors, the approximation is $\tilde U = \sum_{i=1}^{\ranksym{\tilde U}} \tilde v_i \tilde w_i^T$ with $\tilde v_i = d_i \hat v_i$ and the truncation rank $\ranksym{\tilde U}$ is determined by the condition
\begin{align} \label{eq:trunc_criteria}
\sqrt{d^2_{\ranksym{\tilde U}+1} + \cdots + d^2_n } \leq \epsilon_{\text{trunc}} \sqrt{ d_1^2 + \cdots + d_n^2  }.
\end{align}

\subsection{Low-rank GMRES method}
We describe the low-rank GMRES method (lrGMRES) with a generic linear system $Ax=b$. The method follows the standard Arnoldi iteration used by GMRES \cite{saad1986gmres}: construct a set of basis vectors $\{v_i\}_{i=1}^{m_{\text{gm}}}$ by applying the linear operator $A$ to basis vectors, i.e., $w_j = A v_j$ for $j=1,\ldots,m_{\text{gm}}$, and orthogonalizing the resulting vector $w_j$ with respect to previously generated basis vectors $\{v_i\}_{i=1}^{j-1}$. In the low-rank GMRES method, iterates, basis vectors $\{v_i\}$ and intermediate quantities $\{w_i\}$ are represented in terms of the factors of their matricized representations (so that $X$ in \eqref{eq:tensor_rep} would be represented using $Y$ and $Z$ without  explicit construction of $X$), and matrix operations such as matrix-vector products are performed as described in section \ref{sec:matrix_op}. As pointed out in section \ref{sec:matrix_op}, these matrix operations typically tend to increase the rank of the resulting quantity, and this is resolved by interleaving the truncation operator $\mathcal T$ with the matrix operations. The low-rank GMRES method computes a new iterate by solving 
\begin{align}\label{eq:lrgmres_ls}
\min_{\beta \in \mathbb{R}^{m_{\text{gm}}}} \| b - A(x_0 + V_{m_{\text{gm}}} \bar \beta) \|_2,
\end{align}
and constructing a new iterate $x_1 = x_0  + V_{m_{\text{gm}}} \bar \beta$ where $x_0$ is an initial guess. Due to truncation, the basis vectors $\{ v_i\}$ are not orthogonal and $\text{span}(V_{m_{\text{gm}}})$, where $V_{m_{\text{gm}}} = [v_1 \ldots v_{m_{\text{gm}}}]$, is not a Krylov subspace, so that \eqref{eq:lrgmres_ls} 
%the least-squares problem 
must be solved explicitly rather than exploiting Hessenberg structure as in standard GMRES.  Algorithm \ref{alg:lrp_alg} summarizes the lrGMRES. 
We will use this method to solve the linear system of \eqref{eq:jacobian_system}.

\begin{algorithm}[!h]
\caption{Restarted low-rank GMRES method in tensor format \cite{ballani2013projection}}
\label{alg:lrp_alg}
\begin{algorithmic}[1]
\STATE{set the initial solution $\subgmres{\bar {u}^0}$}
\FOR{ $k$ = $0,\,1,\,\dots$} 
\STATE{ $\resiter{}^k := f - A\subgmres{\bar{u}}^k$}
\IF{$\| \resiter{}^k \|_2 / \| f \|_2 < \epsilon_{\text{gmres}}$ \OR $\| \resiter{}^k \|_2 \geq \| \resiter{}^{k-1} \|_2 $}
\STATE{ return $\subgmres{\bar{u}}^k$}
\ENDIF
\STATE{ $\bar {v}_1 := \mathcal T_{\epsilon_{\text{trunc}}}(\resiter{}^k)$} \label{ln:trunc_corr1}
\STATE{ $v_1 := {\bar{v}_1}/{\| \bar{v}_1 \|_2}$}
\FOR{$j = 1,\,\dots,\,m_{\text{gm}}$}
\STATE{ $w_j := Av_j$ \label{lrp-11}}
\STATE{ solve $(V_j^T V_j)\vecrep{\alpha} = V_j^Tw_j$ where $V_j = [v_1,\ldots,v_j]$ }
\STATE{ $\bar{v}_{j+1} :=\mathcal T_{\epsilon_{\text{trunc}}}\left( w_j - \sum_{i=1}^{j} \alpha_i v_i \right) $ } \label{ln:trunc_corr2}
\STATE{ $v_{j+1} := \bar {v}_{j+1} / \| \bar{v}_{j+1}\|_2$ }
\ENDFOR
\STATE{ solve $(W_{m_{\text{gm}}}^T AV_{m_{\text{gm}}}) \vecrep{\beta} = W_{m_{\text{gm}}}^T {r}^k_{\text{gm}}$ where $W_j = [w_1,\ldots,w_j]$}
\STATE{ $\subgmres{\bar{u}}^{k+1} :=\mathcal T_{\epsilon_{\text{trunc}}} (\subgmres{\bar{u}}^k+ V_{m_{\text{gm}}} \bar \beta)$} \label{ln:trunc_corr3}
\ENDFOR
\end{algorithmic}
\end{algorithm}

\subsection{Preconditioning}
We also use preconditioning to speed convergence of the low-rank GMRES method. For this, we consider a right-preconditioned system 
\begin{equation*}
J^n (M^n)^{-1} \tilde u^n = \bar r^n,
\end{equation*}
where $M^n$ is the preconditioner and $M^n \bar u^n = \tilde u^n$ such that $J^n \bar u^n =\bar r^n$. We consider an approximate mean-based preconditioner \cite{powell2009block}, which is derived from the matrix $G_1 \otimes \mathcal F_1$ associated with the mean $\nu_0$ of the random viscosity \eqref{eq:rf_visc},
\begin{equation}\label{eq:precond_matrix}
M^n = G_1 \otimes \begin{bmatrix}M_A^n & B^T \\ 0 & - M^n_s \end{bmatrix},
\end{equation}
where
\begin{alignat*}{2}
M_A^n &= \begin{bmatrix} A_1^{xx} + N_1^{n} & 0 \\ 0& A^{yy}_1+N_1^n \end{bmatrix}, \quad &&(\text{Picard iteration}), \\M_A^n &= \begin{bmatrix} A_1^{xx} + N_1^{n} + W_1^{xx,n} & 0 \\ 0& A^{yy}_1+N_1^n + W_1^{yy,n}\end{bmatrix}, \quad &&(\text{Newton iteration}).
\end{alignat*}
{For approximating the action of the inverse, $(M_s^n)^{-1}$, we choose the boundary-adjusted least-squares commutator (LSC) preconditioning scheme \cite{elman2014finite}, 
\begin{align*}
M_s^n = B F_1^{-1}B^T \approx (B H^{-1} B^T) (B M_\ast^{-1}  F_1 H^{-1}B^T)^{-1} (BM_\ast^{-1}B^T),
\end{align*}
where $M_\ast$ is the diagonal of the velocity mass matrix and $H=D^{-1/2}M_\ast D^{-1/2}$, where $D$ is a diagonal scaling matrix deemphasizing contributions near the boundary.} During the low-rank GMRES iteration, the action of the inverse of the preconditioner \eqref{eq:precond_matrix} can be applied to a vector in a manner analogous to \eqref{eq:mvp_x}--\eqref{eq:mvp_p}.

%% file: adaptive_tolerances.tex
\section{Inexact nonlinear iteration}\label{sec:adap_tol}
As outlined in Algorithm \ref{alg:abs_proc},  we use the hybrid approach, employing a few steps of Picard iteration followed by Newton iteration,  and the linear systems are solved using lrGMRES (Algorithm \ref{alg:lrp_alg}). We extend the hybrid approach to an inexact variant based on an inexact Newton algorithm, in which the accuracy of the approximate linear system solution is tied to the accuracy of the nonlinear iterate (see e.g., \cite{kelley1995iterative} and references therein). That is, when the nonlinear iterate is far from the solution, the linear systems may not have to be solved accurately. Thus, a sequence of iterates  $\bar u^{n+1} \defeq \bar u^{n} + \delta \bar u^{n}$ is computed where $\delta \bar u^{n}$ satisfies
\begin{equation*}
\| J^n_N \delta \bar u^{n} +  \bar r^n \|_2  \leq \epsilon_{\text{gmres}}^n \| \bar r^n \|_2, \quad (J_P \text{ for Picard iteration}),
\end{equation*}
where the lrGMRES stopping tolerance ($\epsilon_{\text{gmres}}^n$ of Algorithm \ref{alg:lrp_alg}) is given by
\begin{equation}\label{eq:gmres_tol}
\epsilon_{\text{gmres}}^n \defeq  \rho_{\text{gmres}} \| \bar r^{n} \|_2,
\end{equation}
where  $0 < \rho_{\text{gmres}} \leq 1$. With this strategy, the Jacobian system is solved with increased accuracy as the error becomes smaller, leading to savings in the average cost per step and, as we will show, with no degradation in the asymptotic convergence rate of the nonlinear iteration. 

In addition, in Algorithms \ref{alg:abs_proc} and \ref{alg:lrp_alg}, the truncation operator $\mathcal T_{\epsilon_{\text{trunc}}}$ is used for the low-rank approximation of the nonlinear iterate (i.e., truncating $\bar u^x$, $\bar u^y$, and $\bar p$ at lines \ref{ln:trunc1} and \ref{ln:trunc2} in Algorithm \ref{alg:abs_proc})  and updates (i.e., truncating $\delta  \bar u^x$, $\delta  \bar u^y$, and $\delta  \bar p$ at lines \ref{ln:trunc_corr1}, \ref{ln:trunc_corr2}, and \ref{ln:trunc_corr3} in Algorithm \ref{alg:lrp_alg}). As the lrGMRES stopping tolerance is adaptively determined by the criterion \eqref{eq:gmres_tol}, we also choose the value of the truncation tolerances $\truncsol{}$ and $\trunccorr{}^n$, adaptively. For truncating the nonlinear iterate, the truncation tolerance for the iterate $\{\truncsol{}^n\}$ is chosen based on the nonlinear iteration stopping tolerance,
\begin{equation*}
\truncsol{} \defeq \rho_{\text{nl}} \truncnl{},
\end{equation*}
where $0<\rho_{\text{nl}}\leq 1$. For truncating the updates (or corrections), the truncation tolerance for the correction $\{\trunccorr{}^n\}$ is adaptively chosen based on the stopping tolerance of the linear solver, 
\begin{align*}
\trunccorr{}^n &\defeq \rho_{\text{trunc,P}} \truncgmres{}^n, \quad (\text{for the $n$th Picard step}),\\
\trunccorr{}^n &\defeq \rho_{\text{trunc,N}} \truncgmres{}^n, \quad (\text{for the $n$th Newton step}),
\end{align*}
where $0 < \rho_{\text{trunc,P}}, \rho_{\text{trunc,N}} \leq 1$. Thus, for computing $n$th update $\delta \bar u^n$, we set $\epsilon_{\text{trunc}} = \trunccorr{}^n$ in Algorithm \ref{alg:lrp_alg}.

\begin{algorithm}[!h]
\caption{Inexact nonlinear iteration with adaptive tolerances}
\label{alg:adaptive_proc}
\begin{algorithmic}[1]
\STATE set $\truncsol{} \defeq \rho_{\text{nl}} \truncnl{} $ 
\STATE compute an approximate solution of $A_{\text{st}} \bar u_{\text{st}} = b_{\text{st}}$  using Algorithm \ref{alg:lrp_alg} \label{alg:stokes_inexact}
\STATE{set an initial guess for the Navier--Stokes problem $\bar u^0 \defeq \bar u_{\text{st}}$}\\
\FOR[Picard iteration]{ $k = 0,\ldots,m_p-1$} \label{alg:picard_start}
\STATE set $\truncgmres{}^k =  \rho_{\text{gmres}} \|\bar r^k\|_2$, and $\trunccorr{}^k = \rho_{\text{trunc,P}} \|\bar r^k\|_2 $ 
\STATE solve $J_{\text{P}}^k\, \delta \bar u^k = -\bar r^k$ using Algorithm \ref{alg:lrp_alg} 
\STATE update $\bar u^{k+1} \defeq \mathcal T_{\truncsol{}} (\bar u^k + \delta \bar u^k)$
\ENDFOR\label{alg:picard_end} \\
\WHILE[Newton iteration]{$\| \bar r^k  \|_2 > \truncnl{} \| \bar r^0\|_2$ } \label{alg:newton_start}
\STATE set $\truncgmres{}^k =  \rho_{\text{gmres}} \|\bar r^k\|_2$, and $\trunccorr{}^k = \rho_{\text{trunc,N}} \|\bar r^k\|_2 $ 
\STATE solve $J_{\text{N}}^k\, \delta \bar u^k = -\bar r^k$ using Algorithm \ref{alg:lrp_alg}
\STATE update $\bar u^{k+1} \defeq \mathcal T_{\truncsol{}} (\bar u^k + \delta \bar u^k)$
\ENDWHILE \label{alg:newton_end}
\end{algorithmic}
\end{algorithm}

%% file: numerical_results2.tex
\section{Numerical results}\label{sec:results}
\begin{figure}[!htb]
	\centering
	\label{fig:fem_domain}
	\includegraphics[angle=0, scale=.5]{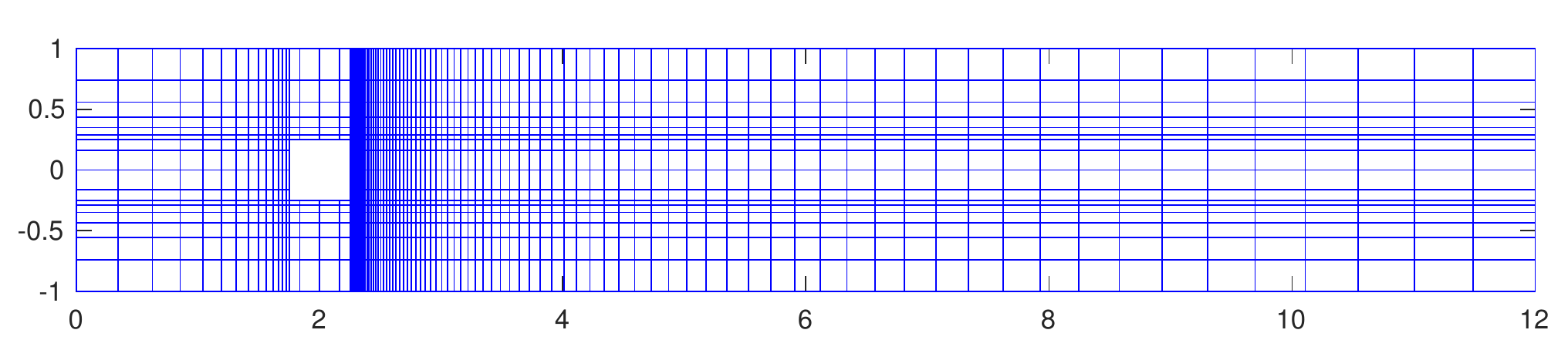}
	\caption{Spatial domain and finite element discretization.}
\end{figure}

In this section, we present the results of numerical experiments on a model problem, flow around a square obstacle in a channel, for which the details are depicted in Figure \ref{fig:fem_domain}. The domain has length 12 and height 2, and it contains a square obstacle centered at (2,0) with sides of length .25. %Figure \ref{fig:fem_domain} depicts the spatial domain and the finite element discretization. 

 For the numerical experiments, we define the random viscosity \eqref{eq:rf_visc} using the Karhunen-Lo\`eve (KL) expansion \cite{loeve1978probability},
\begin{align}\label{eq:visc_KL}
\nu(x,\xi) = \nu_0 + \sigma_\nu \sum_{i=1}^{n_\nu} \sqrt{\lambda_i} \nu_i(x) \xi_i,
\end{align}
where $\nu_0$ and $\sigma^2_\nu$ are the mean and the variance of the viscosity of $\nu(x,\xi)$, and $\{(\lambda_i,\nu_i(x))\}_{i=1}^{n_\nu}$ are eigenpairs of the eigenvalue problem associated with the covariance kernel $C(x,y)$ of the random field. We consider two types of covariance kernel: absolute difference exponential (AE) and squared difference exponential (SE), which are defined via
\begin{align}\label{eq:cov_AE_SE}
C^{\text{AE}}(x,y) = \exp\left(- \sum_{i=1}^2 \frac{ | x_i - y_i|}{l_i} \right), \quad 
C^{\text{SE}}(x,y) = \exp\left(- \sum_{i=1}^2 \frac{(x_i - y_i)^2}{l_i^2}\right),
\end{align}
where $x=(x_1,x_2)$ and $y=(y_1,y_2)$ are points in the spatial domain, and $l_1$, $l_2$ are correlation lengths. We assume that the random variables $\{\xi_i\}_{i=1}^{n_\nu}$ are independent and identically distributed and that $\xi_i$ (for $i=1,\ldots,n_\nu$) follows a uniform distribution over $[-1, 1]$. For the mean of the viscosity, we consider several choices, $\nu_0 = \{\frac{1}{50}, \frac{1}{100}, \frac{1}{150}\}$, which corresponds to $\text{Re}_0 = \{100, 200, 300\}$. In all experiments, we use a finite-term KL-expansion with $n_\nu = 5$. For constructing the finite-dimensional approximation space $S = \text{span}(\{\psi_i(\xi)\}_{i=1}^{n_\xi})$ in the parameter domain, we use orthogonal polynomials $\{\psi_i(\xi)\}_{i=1}^{n_\xi}$ of total degree 3, which results in $n_\xi=56$. The orthogonal polynomials associated with uniform random variables are Legendre polynomials, $\psi_i(\xi) = \prod_{j=1}^{n_\nu} \ell_{d_j(i)}(\xi_j)$ where $d(i) = (d_1(i),\ldots,d_{n_\nu}(i))$ is a multi-index consisting of non-negative integers and $\ell_{d_j(i)}$ is the $d_j(i)$th order Legendre polynomial of $\xi_j$. For the spatial discretization, Taylor--Hood elements are used on a stretched grid, which results in \{6320, 6320, 1640\} degrees of freedom in $\{\vec u^x, \vec u^y, p\}$, respectively (i.e., $n_u =  6320$ and $n_p=1640$.) The implementation is based on the Incompressible Flow and Iterative Solver Software (IFISS) package \cite{elman2014ifiss, ifiss}. %The grid is also shown in Figure \ref{fig:fem_domain}.

\subsection{Low-rank inexact nonlinear iteration}
In this section, we compare the results obtained from the low-rank inexact nonlinear iteration with those obtained from other methods, the exact and the inexact nonlinear iteration with full rank solutions, and the Monte Carlo method. Default parameter settings are listed in Table~\ref{tab:ex_params}, where the truncation tolerances only apply to  the low-rank method.  Unless otherwise specified, the linear system is solved using a restarted version of low-rank GMRES,  lrGMRES(20).%, which generates 20 basis vectors at each cycle.

\begin{table}[htbp]
\caption{Tolerances and adaptive parameters.}
\begin{center}\footnotesize
\renewcommand{\arraystretch}{1.3}
\label{tab:ex_params}
\begin{tabular}{c| c }
\hline
Nonlinear iteration stopping tolerance & $\truncnl{} = 10^{-5}$\\
GMRES tolerance (Stokes) & $\truncgmres{} = 10^{-4}$\\
GMRES tolerances (Picard and Newton) & $\truncgmres{}^n = \rho_{\text{gmres}} \| \bar r^n \|_2$ ($\rho_{\text{gmres}} = 10^{-.5}$)\\
Truncation tolerance for solutions & $\truncsol{} =\rho_{\text{nl}} \truncnl{}$  ($\rho_{\text{nl}}  = 10^{-1}$)\\
Truncation tolerance for corrections & $\trunccorr{}^n =\rho_{\text{trunc}} \truncgmres{}^n$  ($\rho_{\text{trunc}}  = 10^{-1}$)\\
\hline
\end{tabular}
\end{center}
\end{table}

\begin{figure}[!htb]
	\label{fig:InexactNewton_convergence}
	\centering
	\subfloat[Convergence of the hybrid approach]  {
	\includegraphics[angle=0, scale=0.66]{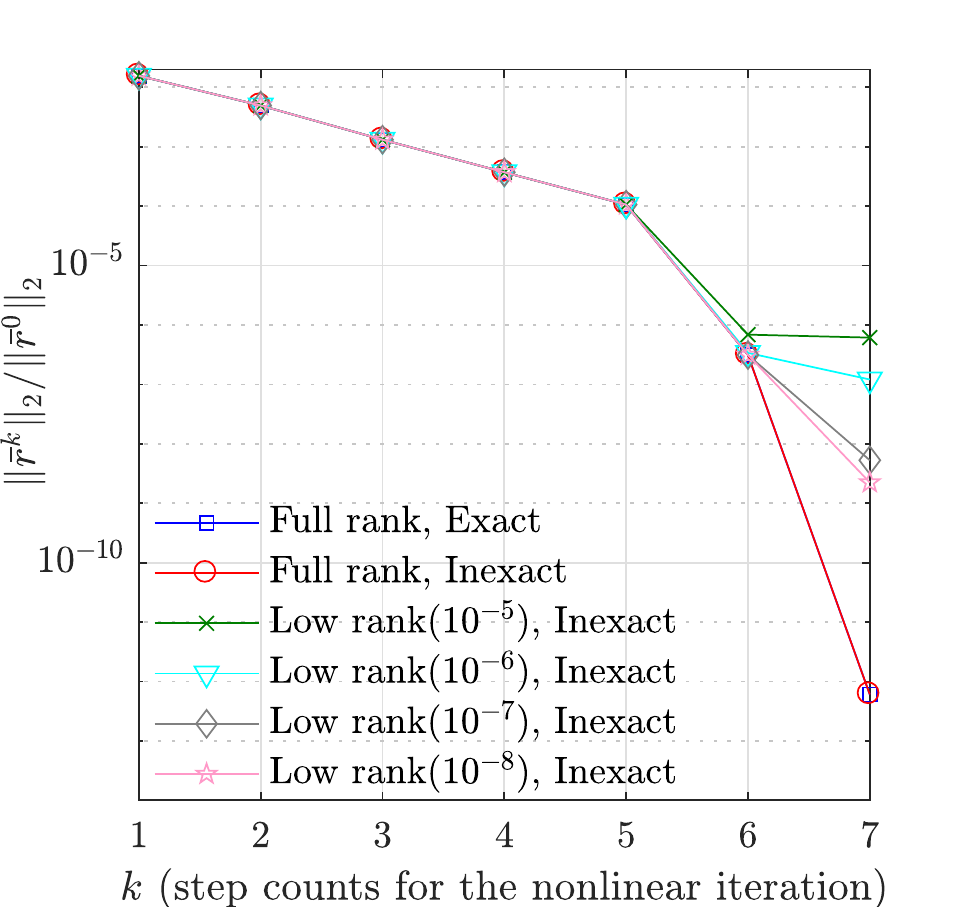}
	\label{fig:conv_newton}
	}
	\subfloat[Convergence of the Picard iteration]  {
	\includegraphics[angle=0, scale=0.66]{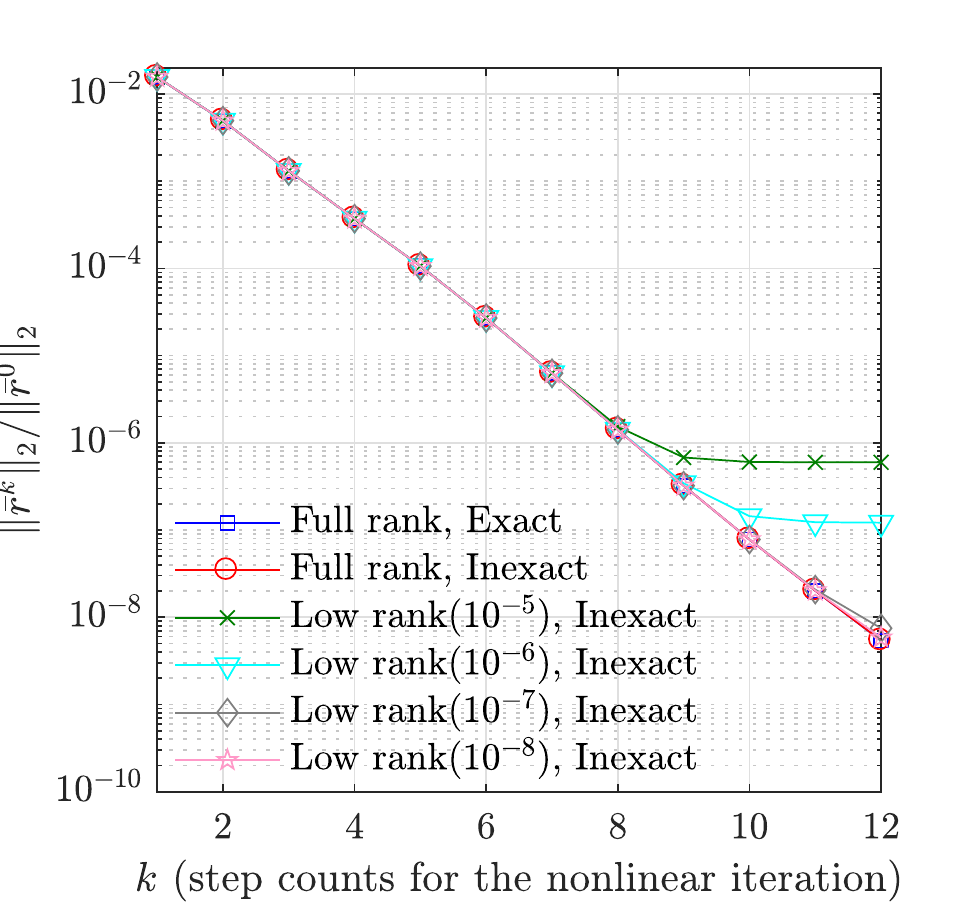}
	\label{fig:conv_picard}
	}
	\caption{Convergence of both exact and inexact nonlinear iterations (full-rank) and the low-rank inexact nonlinear iteration.}
\end{figure}

We first examine the convergence behavior of the inexact nonlinear iteration for a model problem characterized by $\text{Re}_0=100$, $CoV=1\%$, and SE covariance kernel in \eqref{eq:cov_AE_SE} with $l_1=l_2=32$. We compute a full-rank solution using the exact nonlinear iteration ($\truncgmres{}^n = 10^{-12}$ and no truncation) until the nonlinear iterate reaches the nonlinear stopping tolerance, $\truncnl{}=10^{-8}$. Then we compute another full-rank solution using the inexact nonlinear iteration (i.e., adaptive choice of $\truncgmres{}^n$ as shown in Table \ref{tab:ex_params} and no truncation). Lastly, we compute a low-rank approximate solution using the low-rank inexact nonlinear iteration (i.e., adaptive choices of $\truncgmres{}^n$ and $\trunccorr{}^n$ as shown in Table \ref{tab:ex_params} and for varying $\truncsol{}=\{ 10^{-5},10^{-6},10^{-7},10^{-8}\}$). Figure \ref{fig:InexactNewton_convergence} shows the convergence behavior of the three methods. In Figure \ref{fig:conv_newton}, the hybrid approach is used, in which the first step corresponds to the Stokes problem (line \ref{alg:stokes_inexact} of Algorithm~\ref{alg:adaptive_proc}), the 2nd--5th steps correspond to the Picard iteration (line \ref{alg:picard_start}--\ref{alg:picard_end} of Algorithm \ref{alg:adaptive_proc}, and $m_p=4$), and the 6th--7th steps correspond to the Newton iteration (line \ref{alg:newton_start}--\ref{alg:newton_end} of Algorithm \ref{alg:adaptive_proc}). Figure \ref{fig:conv_newton} confirms that the inexact nonlinear iteration is as effective as the exact nonlinear iteration. The low-rank inexact nonlinear iteration behaves similarly up to the 6th nonlinear step but when the truncation tolerances are large $\truncsol{}=\{10^{-5}, 10^{-6}\}$, it fails to produce a nonlinear solution satisfying $\truncnl{}=10^{-8}$. 
%, which is due to the mild truncation tolerances $\truncsol{}=\{10^{-5}, 10^{-6}\}$. 
Similar results can be seen in Figure \ref{fig:conv_picard}, where only the Picard iteration is used. As expected, in that case, the relative residual decreases linearly for all solution methods, but the low-rank inexact nonlinear iteration with the mild truncation tolerances also fails to reach the nonlinear iteration stopping tolerance.

\begin{figure}[!tb]
\hspace{-27mm}
	\label{fig:Re100Cov1Corrlen3200SEfull_stats_aug}
	\includegraphics[angle=0, width=175mm, height=58mm]{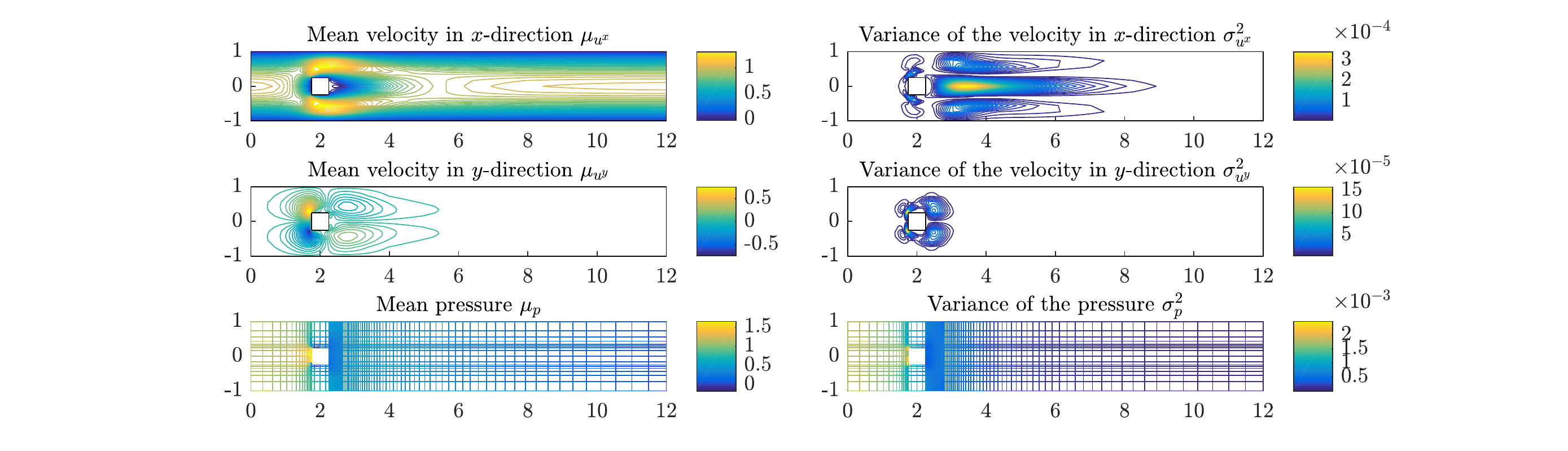}
	\vspace{-8mm}
	\caption{Mean and variances of full-rank velocity solutions $\vec u^x(x,\xi) $, $\vec u^y(x,\xi)$, and pressure solution $p(x,\xi)$ for {\normalfont $\text{Re}_0=100$}, $CoV=1$, and $l_1=l_2=32$.}
\end{figure}

\begin{figure}[!tb]
\hspace{-27mm}
	\label{fig:Re100Cov1Corrlen3200SEdiff_stats_aug}
	\includegraphics[angle=0, width=175mm, height=58mm]{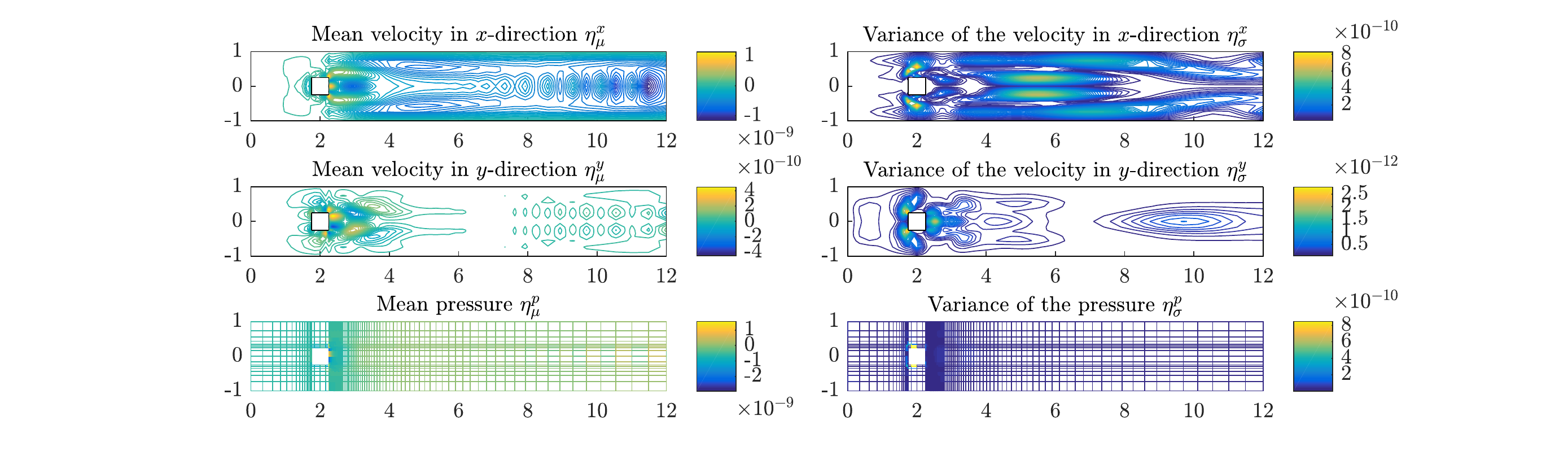}
	\vspace{-8mm}
	\caption{Difference in the means and variances of the full-rank and the low-rank solutions for {\normalfont $\text{Re}_0=100$}, $CoV=1$, and $l_1=l_2=32$.}
\end{figure}

Figure \ref{fig:Re100Cov1Corrlen3200SEfull_stats_aug} shows means and variances of the components of the full-rank solution, given by
\begin{alignat}{5}
% u_0^x &= E[ \vec u^x ], \quad  && u_0^y &&= E[ \vec u^y ], \quad && p_0 &&= E[ p ], \label{eq:mean_gpc}\\
 \mu_{u^x} &= \mathbb E[ \vec u^x ], \quad  && \mu_{u^y} &&= \mathbb E[ \vec u^y ], \quad && \mu_p &&= \mathbb E[ p ], \label{eq:mean_gpc}\\
\sigma^2_{u^x} &= \mathbb E[ (\vec u^x  -  \mu_{u^x})^2], \quad &&\sigma^2_{u^y} &&= \mathbb E[ (\vec u^y  -  \mu_{u^y})^2],\quad &&\sigma^2_{p} &&= \mathbb E[ (p  -  \mu_{p})^2]\label{eq:var_gpc}.
\end{alignat}  
These quantities are easily computed by exploiting the orthogonality of basis functions in the gPC expansion. Figure \ref{fig:Re100Cov1Corrlen3200SEdiff_stats_aug} shows the differences in the means and variances of the solutions computed using the full-rank and the low-rank inexact nonlinear iteration. Let us denote the full-rank and low-rank horizontal velocity solutions by $u^{x,\text{full}}$ and $u^{x,\text{lr}}$, with analogous notation for the vertical velocity and the pressure. Thus, the differences in the means and the variances are 
\begin{alignat*}{5}
% \eta_0^x &=  u_0^{x,\text{full}} -  u_0^{x,\text{lr}}, \quad  &&\eta_0^y &&=  u_0^{y,\text{full}} -  u_0^{y,\text{lr}}, \quad && \eta^p_0 &&= p_0^{\text{full}} - p_0^{\text{lr}}, \\
\eta_\mu^x &=  \mu_{u^{x,\text{full}}} -  \mu_{u^{x,\text{lr}}} , \quad  &&\eta_\mu^y &&=  \mu_{u^{y,\text{full}}} -  \mu_{u^{y,\text{lr}}}, \quad && \eta_\mu^p &&= \mu_{p^{\text{full}}} -  \mu_{p^{\text{lr}}}, \\
\eta_\sigma^x &= \sigma^2_{u^{x,\text{full}}} - \sigma^2_{u^{x,\text{lr}}}, \quad && \eta_\sigma^y &&= \sigma^2_{u^{y,\text{full}}} - \sigma^2_{u^{y,\text{lr}}},\quad &&\eta_\sigma^p &&= \sigma^2_{p^{\text{full}}} - \sigma^2_{p^{\text{lr}}}.
\end{alignat*}  
Figure \ref{fig:Re100Cov1Corrlen3200SEdiff_stats_aug} shows these differences, normalized by graph norms $\| \nabla \vec  \mu_{u^{\text{full}}} \| + \| \mu_{p^{\text{full}}} \| $ for the means and $\| \nabla \vec \sigma_{u^{\text{full}}}^2 \| + \| \sigma^2_{p^{\text{full}}} \|$ for the variances, where  $\| \nabla \vec u \|  = (\int_D \nabla \vec u :  \nabla \vec u \,dx)^{\frac{1}{2}} $ and $\| p \| =  (\int_D p^2 dx)^{\frac{1}{2}}$.
%normalized versions of these differences, in which each difference is divided by the graph norm of the full-rank solution, $\| \nabla \vec u \| + \| p \| = (\int_D \nabla \vec u :  \nabla \vec u \,dx)^{\frac{1}{2}} + (\int_D p^2 dx)^{\frac{1}{2}}$. 
Figure \ref{fig:Re100Cov1Corrlen3200SEdiff_stats_aug} shows that the normalized differences in the mean and the variance  are of order $10^{-9}\sim10^{-10}$ and $10^{-10}\sim10^{-12}$, respectively, i.e., the errors in low-rank solutions are considerably smaller than the magnitude of the truncation tolerances $\truncsol{}, \trunccorr{}$ (see Table \ref{tab:ex_params}).

\subsection{Characteristics of the Galerkin solution} 
In this section, we examine various properties of the Galerkin solutions, with emphasis on comparison of the low-rank and full-rank versions of these solutions and development of an enhanced understanding of the relation between the Galerkin solution and the polynomial chaos basis. We use the same experimental setting studied above (SE covariance kernel, $l_1=l_2=32$, $\text{Re}_0=100$ and $CoV=1\%$).

\begin{figure}[!tb]
	\centering
	\label{fig:qoi_pdf1}
	\subfloat[$\vec u^x$]  {
	\includegraphics[angle=0, scale=0.60]{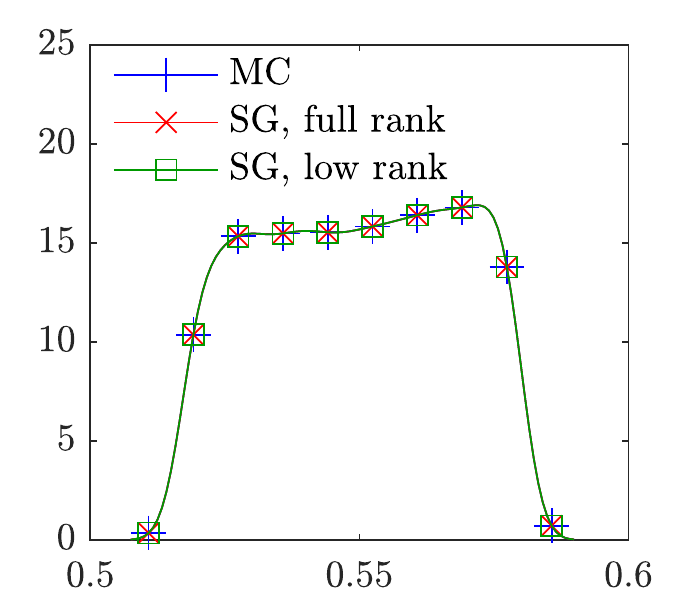}
	\label{fig:qoi_ux1}
	}
	\subfloat[$\vec u^y$]  {
	\includegraphics[angle=0, scale=0.60]{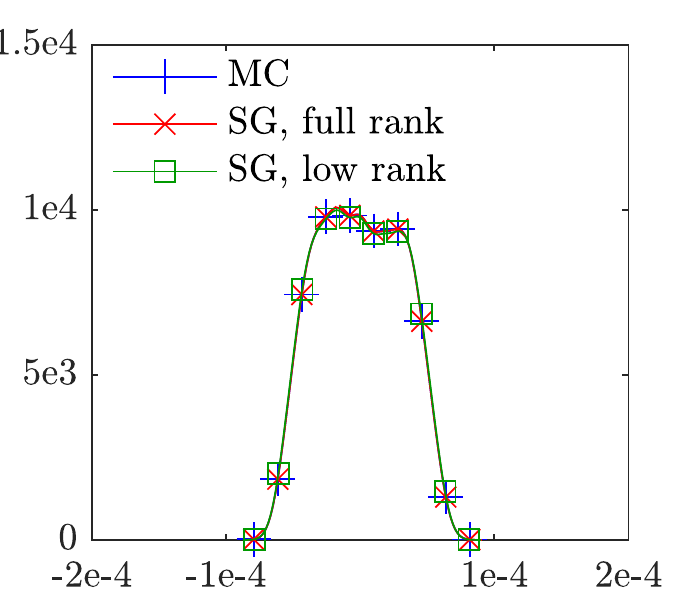}
	\label{fig:qoi_uy1}
	}
	\subfloat[$p$]  {
	\includegraphics[angle=0, scale=0.60]{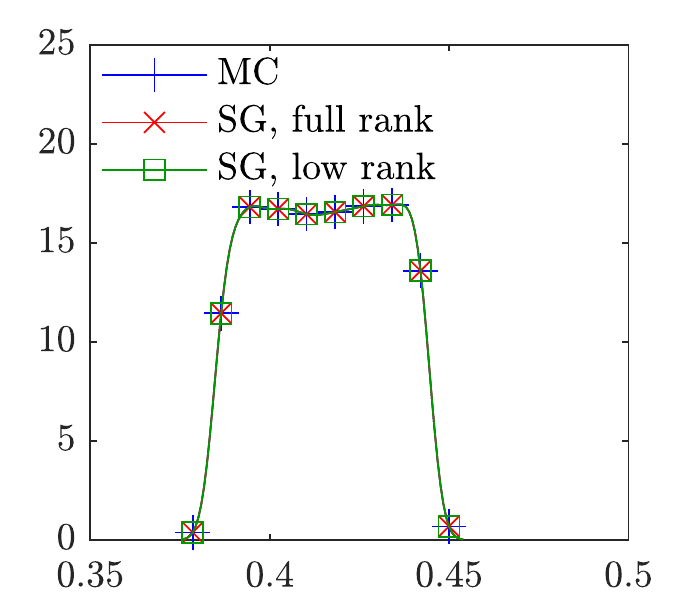}
	\label{fig:qoi_p1}
	}
	\caption{Estimated pdfs of the velocities $\vec u^x$, $\vec u^y$, and the pressure $p$ at the point (3.6436, 0).}
\end{figure}

We begin by comparing the Galerkin solution with one obtained using Monte Carlo methods. In particular, we estimate a probability density function (pdf) of the velocity solutions ($\vec u^x(x,\xi)$, $\vec u^y(x,\xi)$) and the pressure solution ($p(x,\xi)$) at a specific point on the spatial domain $D$. In the Monte Carlo method, we solve $n_{\text{MC}} = 25000$ deterministic systems, $\nlres{\vec u,p,\vec v,q;\xi^{(k)}}=0$ associated with $n_{\text{MC}}$ realizations $\{\xi^{(k)}\}_{k=1}^{n_{\text{MC}}}$ in the parameter space. Using the \textsc{Matlab} function {\tt ksdensity}, the pdfs of ($\vec u^x(x,\xi)$, $\vec u^y(x,\xi)$, $p(x,\xi)$) are estimated at the spatial point with coordinates (3.6436, 0), where the variance of $\vec u^x(x,\xi)$ is large (see Figure \ref{fig:Re100Cov1Corrlen3200SEfull_stats_aug}). The results are shown in Figure \ref{fig:qoi_pdf1}. They indicate that the pdf of the Galerkin solution is virtually identical to that of the Monte Carlo solution, and there is essentially no difference between the low-rank and full-rank results.

\begin{figure}[!tb]
	\centering
	\label{fig:Re100Cov1Corrlen3200SEfull_gPC_decay}
	\includegraphics[angle=0, scale=.75]{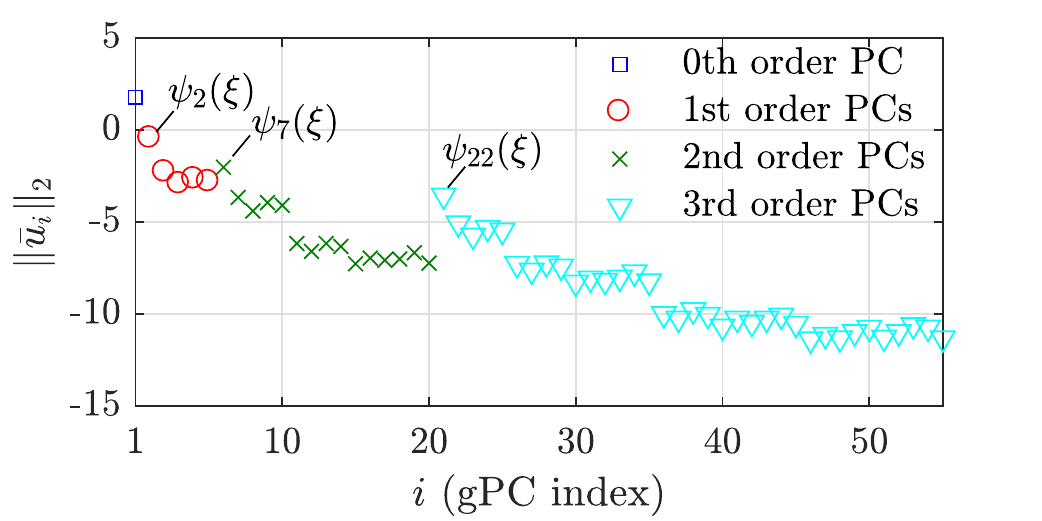}
	\caption{Norms of the gPC coefficients $\|\bar u_i\|_2$ for {\normalfont $\text{Re}_0=100$}, $CoV=1$, and $l_1=l_2=32$.}
\end{figure} 

Next, we explore some characteristics of the Galerkin solution, focusing on the horizontal velocity solution; the observations made here also hold for the other components of the solution.  Given the coefficients of the velocity solution in matricized form, $U^x$, the discrete velocity solution is then given by 
\begin{equation*}
\vec u^x(x,\xi) = \Phi^T(x) U^x \Psi(\xi),
\end{equation*}
where $\Phi(x) = [\phi_1(x),\ldots,\phi_{n_u}(x)]^T$ and $\Psi(\xi) = [\psi_1(\xi),\ldots,\psi_{n_\xi}(\xi)]^T$. Consider in particular the component of this expression corresponding to the jth column of $U^x$,
$$
\left(\sum_{i=1}^{n_u} u_{ij}^x \phi_i(x) \right) \psi_j(\xi)
$$
so that this ($j$th) column $\bar u_j^x = [U^x]_j$ corresponds to the coefficient of the $j$th polynomial basis function $\psi_j$. Figure \ref{fig:Re100Cov1Corrlen3200SEfull_gPC_decay} plots the values of the coefficients $\| \bar u^x_i \|_2$. (This data is computed with $\text{Re}_0=100$, $CoV=1\%$, and SE covariance kernel with $l_1=l_2=32$). Note that the gPC indices $\{j\}$ are in one-to-one correspondence with 
multi-indices $d(j) = (d_1(j), \ldots,d_{n_u}(j))$, where the element of the multi-index indicates the degree of univariate Legendre polynomial. The multi-indices $\{d(i)\}_{i=1}^{n_\xi}$ are ordered in the lexicographical order, for example, the first eight multi-indices are as $d(1) = (0,0,0,0,0),\, d(2) = (1,0,0,0,0),\,  d(3) = (0,1,0,0,0),\,\ldots,\, d(6) = (0,0,0,0,1),\, d(7) = (2,0,0,0,0),$ and $d(8) = (1,1,0,0,0)$. In Figure \ref{fig:Re100Cov1Corrlen3200SEfull_gPC_decay}, the blue square is associated with the zeroth-order gPC component ($d(1)$), the red circles are associated with the first-order gPC components ($\{d(i)\}_{i=2}^{6}$), and so on. Let us focus on three gPC components associated only with $\xi_1$, \{$\psi_2(\xi)=\ell_1(\xi_1),\, \psi_7(\xi)=\ell_2(\xi_1),\, \psi_{22}(\xi)=\ell_3(\xi_1)\}$, where, for $i=2,7,22$, the multi-indices are $d(2) = (1,0,0,0,0)$, $d(7)=(2,0,0,0,0)$, and $d(22)=(3,0,0,0,0)$. The figure shows that the coefficients of gPC components \{$\psi_2(\xi), \psi_7(\xi), \psi_{22}(\xi)\}$ decay more slowly than those of gPC components associated with other random variables $\{\xi_i\}_{i=2}^{n_\nu}$.

\begin{figure}[!tb]
	\subfloat[Plots of coefficients of gPC components 2--7 of $\vec u^x(x,\xi)$]  {
	\hspace{-27mm}
	\label{fig:Re100Cov1Corrlen3200SEfull_gPC}
	\includegraphics[angle=0, width=175mm, height=58mm]{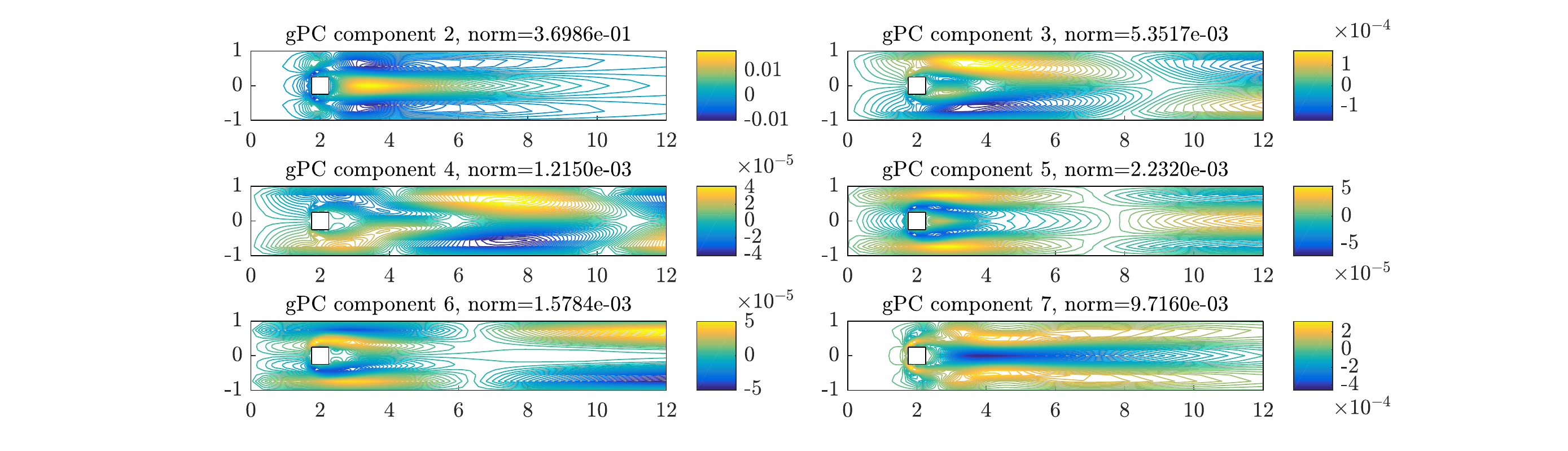}
	\vspace{-8mm}
	}\\
	\subfloat[Plots of coefficients $v_i$ of $\theta_i^x(\xi)$ for $i=2,\ldots,7$]  {
	\hspace{-27mm}
	\label{fig:Re100Cov1Corrlen3200SElr_gPC}
	\includegraphics[angle=0, width=175mm, height=58mm]{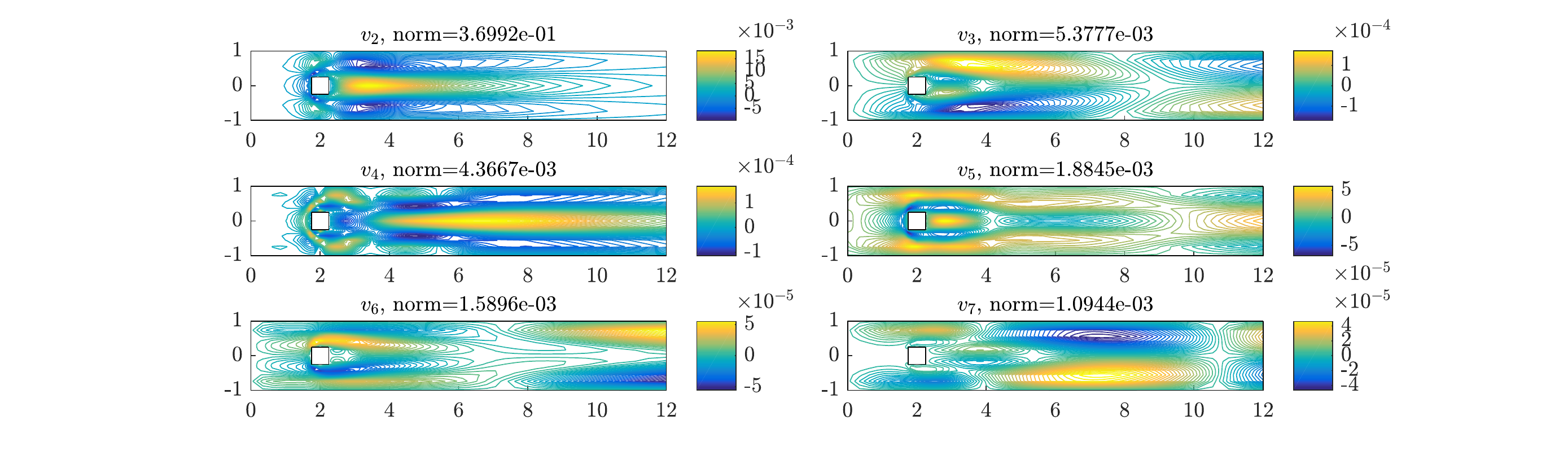}
	\vspace{-8mm}
	}	
	\caption{Plots of coefficients of gPC components 2--7 of $\vec u^x(x,\xi)$ and coefficients $v_i$ of $\theta_i(\xi)$ for $i=2,\ldots,7$ for {\normalfont $\text{Re}_0=100$}, $CoV=1$, and $l_1=l_2=32$.}
\end{figure}

We continue the examination of this data in Figure  \ref{fig:Re100Cov1Corrlen3200SEfull_gPC}, which shows two-dimensional mesh plots of the 2nd through 7th columns of $U^x$. These images show that these coefficients are either symmetric with  respect to the horizontal axis, or reflectionally symmetric (equal in magnitude but of opposite sign), and (as also revealed in Figure \ref{fig:Re100Cov1Corrlen3200SEfull_gPC_decay}), they tend to have smaller values as the index $j$ is increased.

We now look more closely at features of the factors of the low-rank approximate solution and compare these with those of the (unfactored) full-rank solution.  In the low-rank format, the solution is represented using factors $\vec u^x\!(x,\xi)\!=\!(\Phi^T\!(x)V^x)(\Psi^T\!(\xi)W^x)^T\!\!.$ Let us introduce a concise notation of $\vec u^x(x,\xi)\!=\!Z^x_{\ranksym{\bar u^x}}\!(x)^T \Theta^x_{\ranksym{\bar u^x}}\!(\xi) = \sum_{i=1}^{\ranksym{\bar u^x}} \zeta^x_i(x) \theta^x_i(\xi)$ where $Z^x_{\ranksym{\bar u^x}}(x) = [\zeta^x_1(x), \ldots, \zeta^x_{\ranksym{\bar u^x}}(x)]$ and $\Theta^x_{\ranksym{\bar u^x}}(\xi) = [\theta^x_1(\xi), \ldots, \theta^x_{\ranksym{\bar u^x}}(\xi)]$ with $\zeta^x_i(x) = [\Phi^T(x) V^x]_i$ and $\theta^x_i(\xi) = [(\Psi^T(\xi) W^x)]_i$ for $i=1,\ldots,\ranksym{\bar u^x}$. Figure \ref{fig:Re100Cov1Corrlen3200SElr_gPC} shows the coefficients of the $i$th random variable $\theta_i(\xi)$. As opposed to the gPC coefficients of the full-rank solution, the norms of the coefficients of $\{\theta_i(\xi)\}$ decrease monotonically as the index $i$ increases. This is a consequence of the fact that the ordering for $\{\theta_i(\xi)\}$ comes from the singular values of $U^x$. Figure  \ref{fig:Re100Cov1Corrlen3200SElr_gPC} shows the 2nd-7th columns of $V^x$. Figures \ref{fig:Re100Cov1Corrlen3200SEfull_gPC} and \ref{fig:Re100Cov1Corrlen3200SElr_gPC} show that the coefficients $\{v_i\}$ of $\{\theta_i(\xi)\}$ are comparable to the coefficients $\{u_i^x\}$ of the gPC components. Each pair of components in the following parentheses is similar to each other: $(u_2, v_2)$, $(u_3, v_3)$, $(u_7, -v_4)$, $(u_4,-v_7)$, $(u_5,v_5)$, and $(u_6, -v_6)$.

\begin{figure}[!tb]
\centering
	\label{fig:Re100Cov1Corrlen3200SElr_heatmap}
	\includegraphics[angle=0, scale=.8]{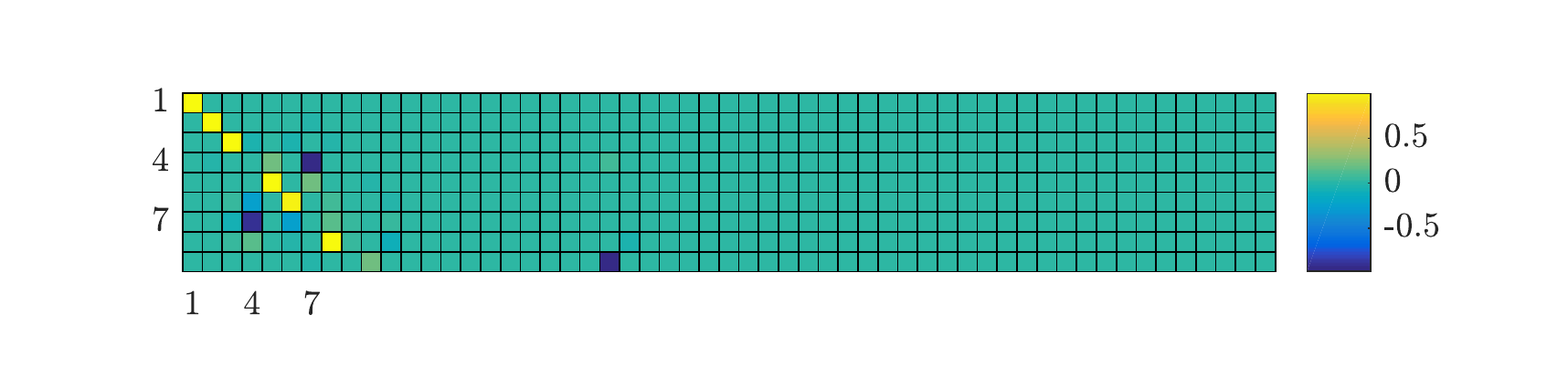}
	\vspace{-10mm}
	\caption{A heat map of $(W^x)^T$.}
\end{figure}

While the columns of $V^x$ show the resemblance to the subset of the columns of $U^x$, $W^x$ tends to act as a permutation matrix. Figure \ref{fig:Re100Cov1Corrlen3200SElr_heatmap} shows a ``heat map'' of $(W^x)^T$, where values of the elements in $W^x$ are represented as colors and the map shows that a very few elements of $W^x_i$ are dominant and a sum of those elements is close to 1. Recall that $\theta^x_i(\xi) =  (W^x_i)^T \Psi(\xi)$. Many dominant elements are located in the diagonal of $W^x$, which results in $\theta^x_i(\xi) \approx \pm \psi_i(\xi)$ (e.g., $i=1,2,3,5,\ldots$). In the case of $W^x_4$, the most dominant element is the 7th element and has a value close to -1, which results in $\theta^x_4(\xi) \approx - \psi_7(\xi)$. As observed in Figure \ref{fig:Re100Cov1Corrlen3200SEfull_gPC_decay}, $\psi_7(\xi)$ has a larger contribution than other gPC components and, in the new solution representation, $\theta^x_4(\xi)$, which consists mainly of $\psi_7(\xi)$, appears earlier in the representation.

\subsection{Computational costs} In this section, we assess the costs of the low-rank inexact nonlinear iteration under various experimental settings: two types of covariance kernels \eqref{eq:cov_AE_SE}, varying $CoV$ \eqref{eq:cov_def}, and varying $\text{Re}_0$. In addition, for various values of these quantities, we investigate the decay of the eigenvalues $\{ \lambda_i\}$ used to define the random viscosity \eqref{eq:visc_KL} and their influence on the rank of solutions. All numerical experiments are performed on an Intel 3.1 GHz i7 CPU, 16 GB RAM using \textsc{Matlab} R2016b and costs are measured in terms of CPU wall time (in seconds). For larger $CoV$ and $\text{Re}_0$, we found the solver to be more effective using the slightly smaller truncation tolerance $\rho_{\text{trunc}} = 10^{-1.5}$ and used this choice for all experiments described below. (Other adaptive tolerances are those shown as in Table \ref{tab:ex_params}.) This change had little impact on results for small $CoV$ and $\text{Re}_0$.

%\REV{To cope with difficulties that may arise from larger $CoV$ or higher $\text{Re}_0$, we set $\rho_{\text{trunc}} = 10^{-1.5}$ for the following experiments. Other adaptive tolerances are those shown as in Table~\ref{tab:ex_params}.}

\begin{figure}[!h]
	\centering
	\includegraphics[angle=0, scale=0.61]{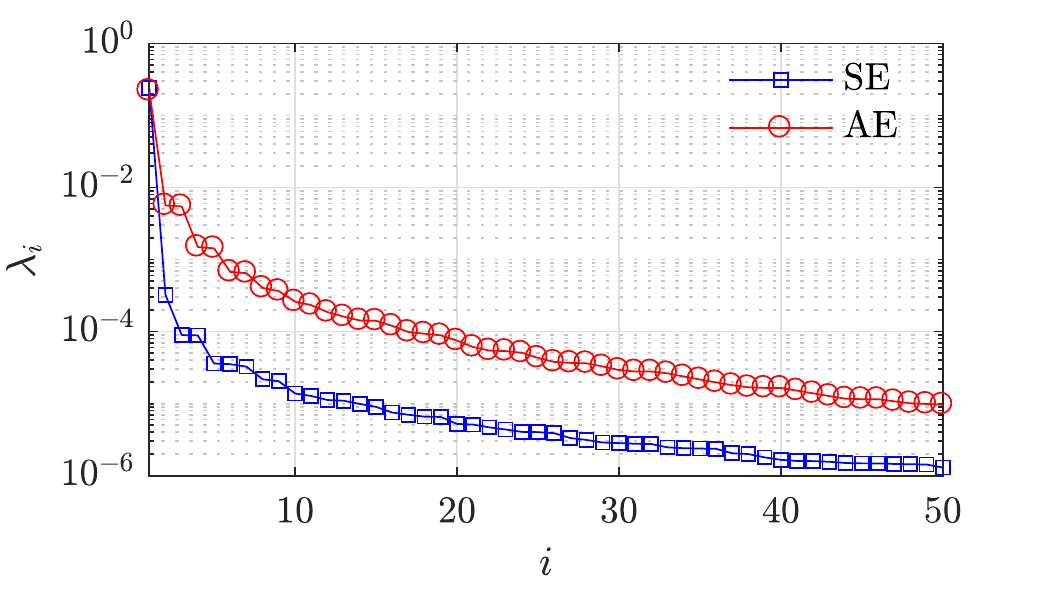}
	\caption{Eigenvalue decay of the AE and the SE covariance kernels.}
	\label{fig:eig_decay_varying_kernel}
\end{figure}

Figure \ref{fig:eig_decay_varying_kernel} shows the 50 largest eigenvalues $\{ \lambda_i\}$ of the eigenvalue problems associated with the SE covariance kernel and the AE covariance kernel \eqref{eq:cov_AE_SE}  with $l_1=l_2=8$, $CoV=1\%$, and $\text{Re}_0=100$. The eigenvalues of the SE covariance kernel decay much more rapidly than those of the AE covariance kernel. Because we choose a fixed number of terms $n_\nu=5$, the random viscosity with the SE covariance kernel retains a smaller variance. 

\begin{figure}[!h]
	\centering
	\subfloat[Computational cost of full-rank computation and low-rank approximation]  {
	\includegraphics[angle=0, scale=0.67]{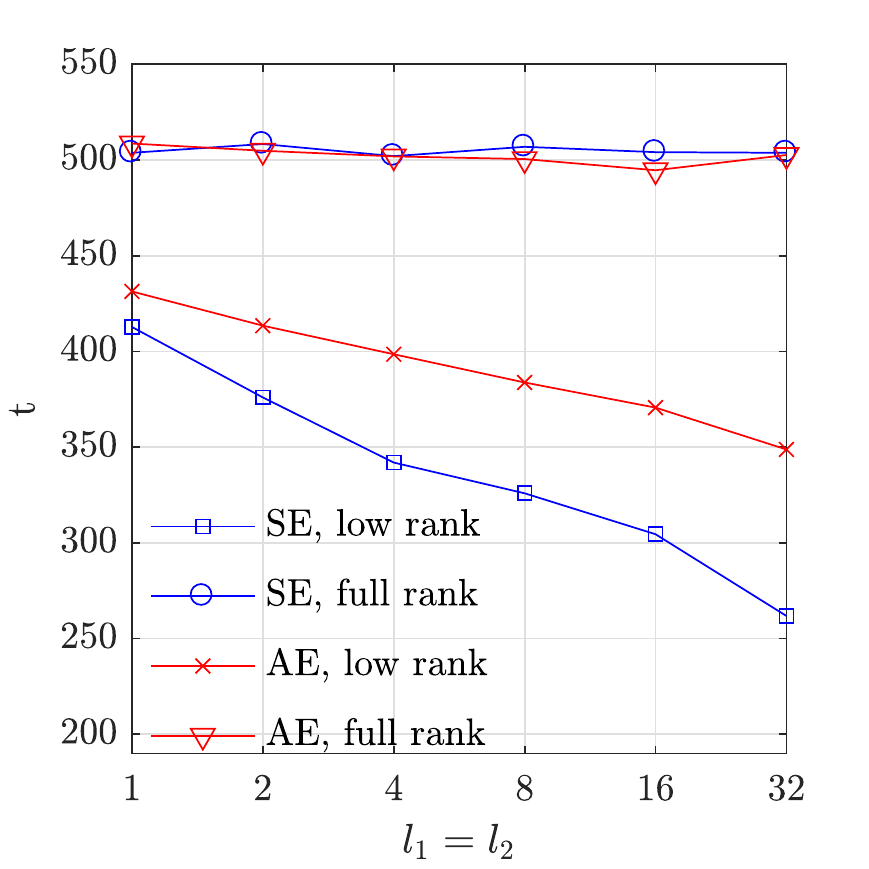}
	\label{fig:cost_AE_SE}
	}\hspace{3mm}	
	\subfloat[Ranks of the low-rank approximate solutions]  {
	\includegraphics[angle=0, scale=0.67]{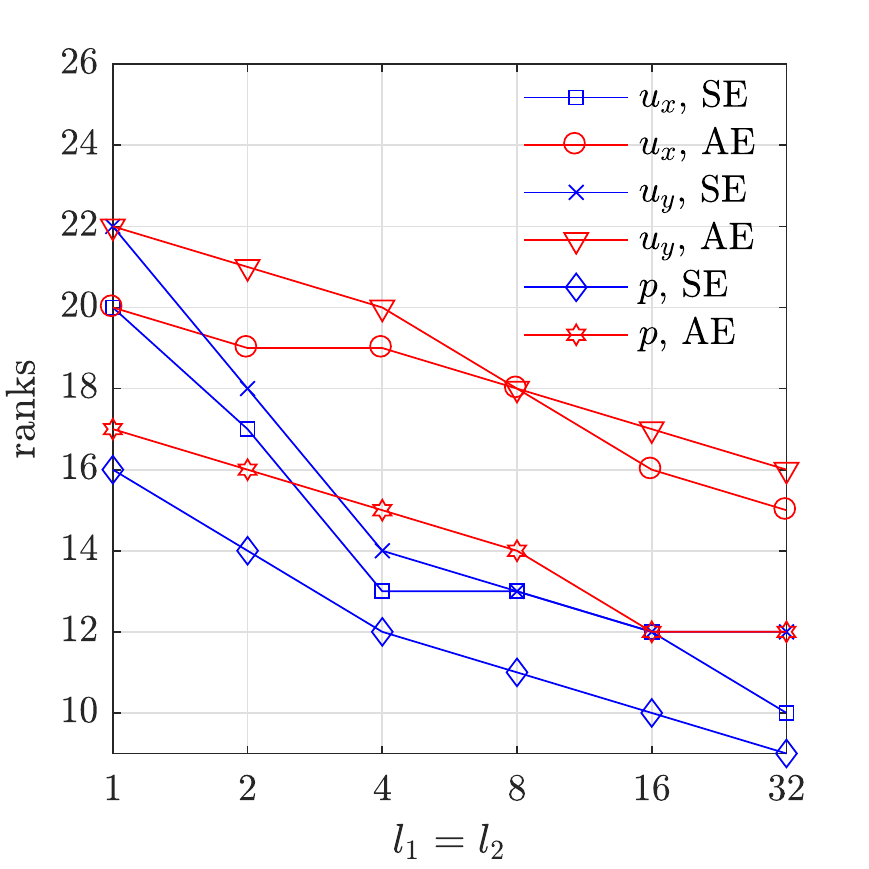}
	\label{fig:rank_AE_SE}
	}
	\caption{Computational costs and ranks for varying correlation lengths with SE and AE covariance kernel.}
\end{figure}

Figure \ref{fig:cost_AE_SE} shows the computational costs (in seconds) needed  for computing the full-rank solutions and the low-rank approximate solutions using the inexact nonlinear iteration for the two covariance kernels and a set of correlation lengths, $l_1 = l_2 = \{1, 2, 4, 8, 16, 32\}$. Figure \ref{fig:rank_AE_SE} shows the ranks of the low-rank approximate solutions that satisfy the nonlinear stopping tolerance $\epsilon_{\truncnl{}}=10^{-5}$.  Again, $\text{Re}_0 = 100$ and $CoV = 1\%$. For this benchmark problem, 4 Picard iterations and 1 Newton iteration are enough to generate a nonlinear iterate satisfying the stopping tolerance $\truncnl{}$. It can be seen from Figure \ref{fig:cost_AE_SE} that in all cases the use of low rank methods reduces computational cost. Moreover, as the correlation length becomes larger, the ranks of the corrections and the nonlinear iterates become smaller. As a result, the low-rank method achieves greater computational savings for the problems with larger correlation length.

Next, we examine the performances of the low-rank approximation method for varying $CoV$, which is defined in \eqref{eq:cov_def}. In this experiment, we fix the value of $\text{Re}_0=100$ and the variance of the random $\sigma_\nu$ is controlled. We consider the SE covariance kernel. 

\begin{figure}[!htb]
	\label{fig:SE_varyingCoV}
	\centering
	\subfloat[Computational cost of full-rank computation and low-rank approximation]  {
	\includegraphics[angle=0, scale=0.60]{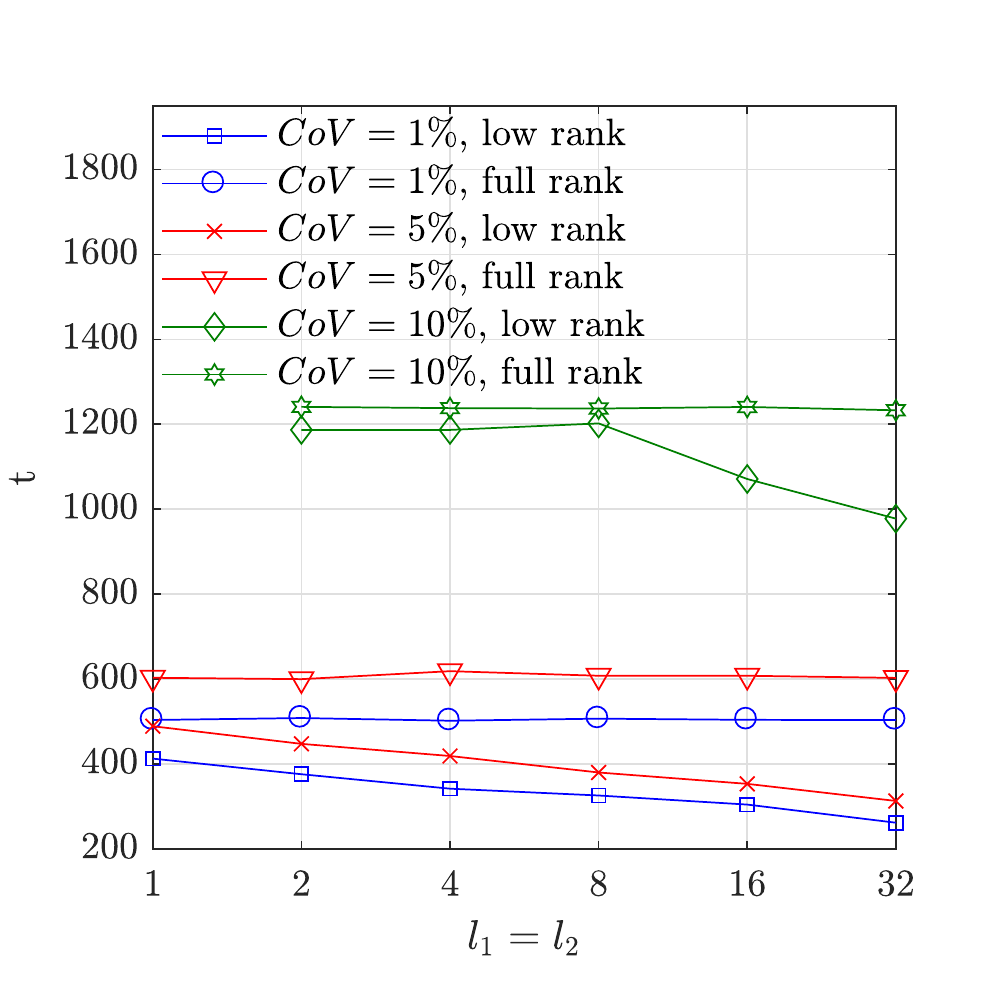}
	\label{fig:cost_SE_varyingCoV}
	}\hspace{3mm}	
	\subfloat[Ranks of the low-rank approximate solutions $u_x$]  {
	\includegraphics[angle=0, scale=0.60]{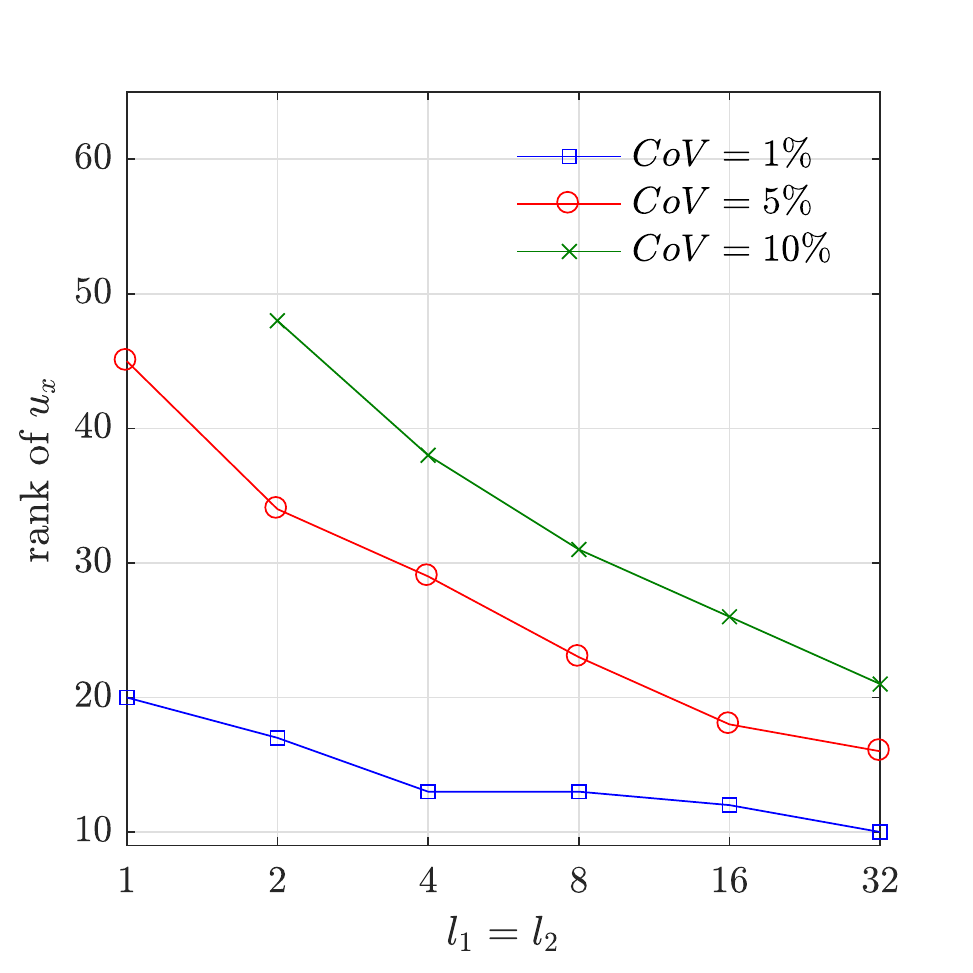}
	\label{fig:rank_SE_varyingCoV_ux}
	}\\
	\subfloat[Ranks of the low-rank approximate solutions $u_y$]  {
	\includegraphics[angle=0, scale=0.60]{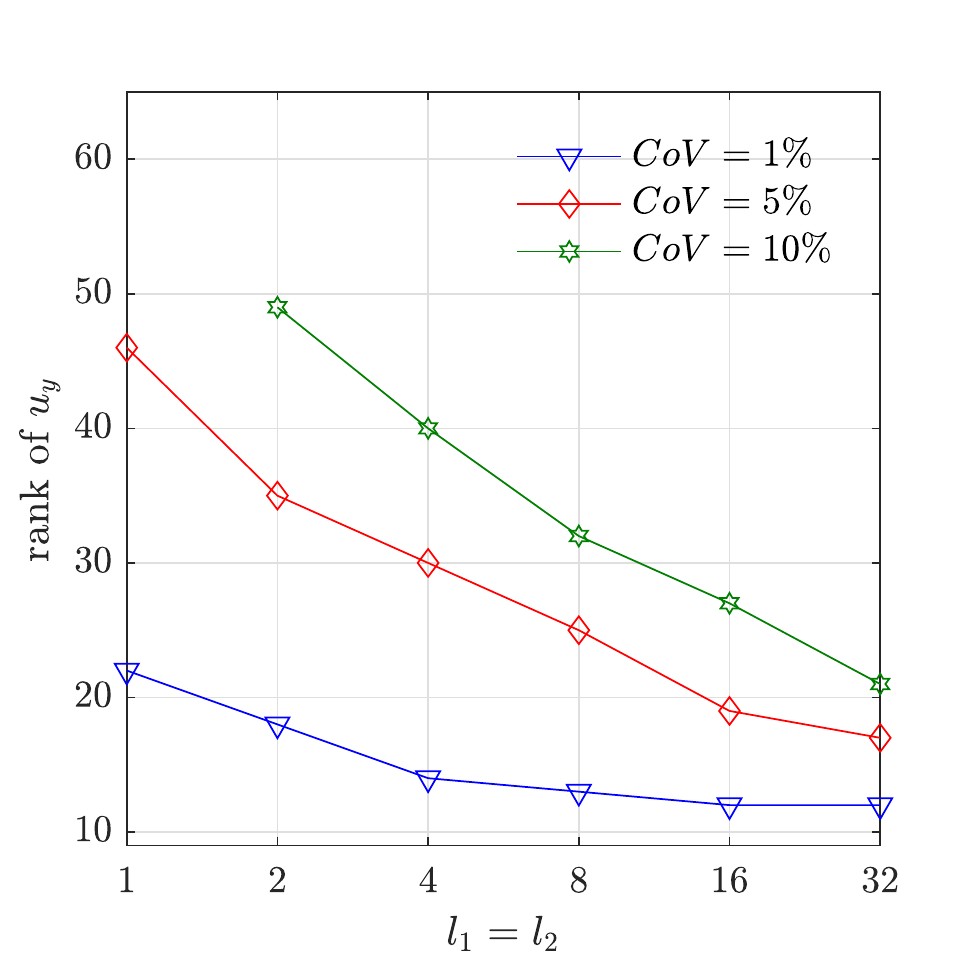}
	\label{fig:rank_SE_varyingCoV_uy}
	}\hspace{3mm}	
	\subfloat[Ranks of the low-rank approximate solutions $p$]  {
	\includegraphics[angle=0, scale=0.60]{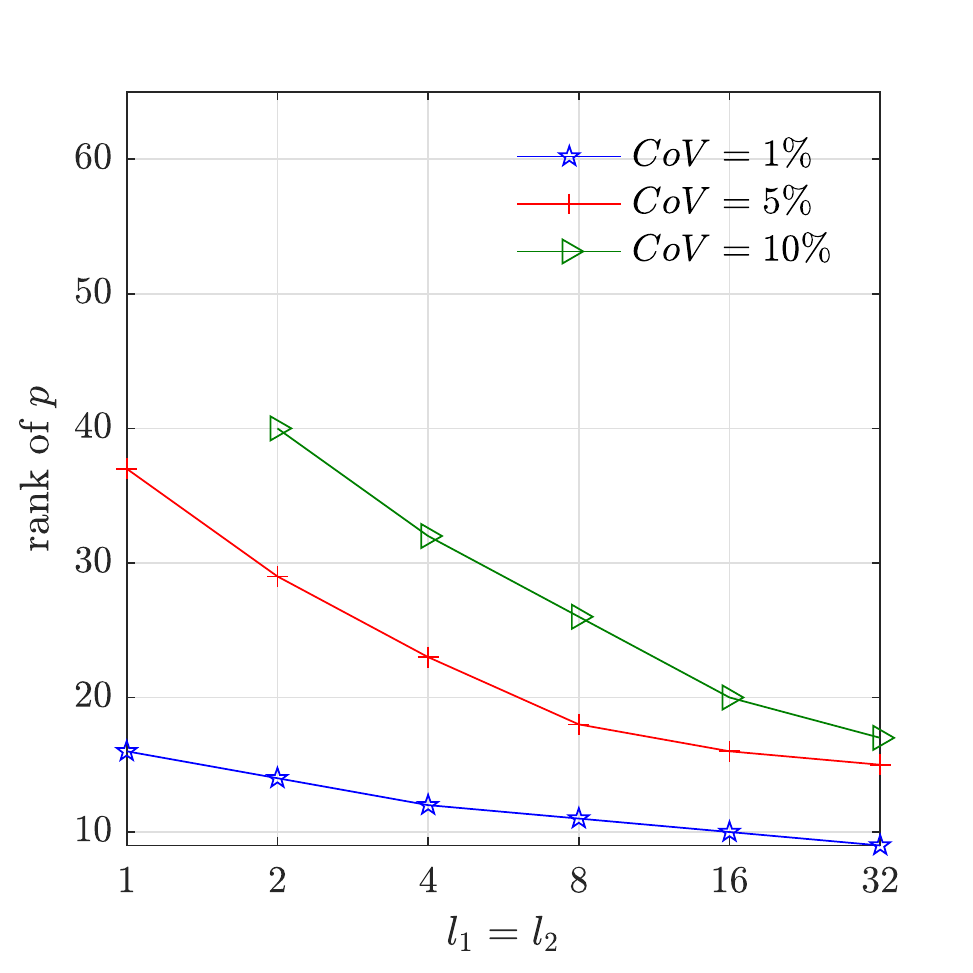}
	\label{fig:rank_SE_varyingCoV_p}
	}
	\caption{Computational costs and ranks for varying correlation lengths and varying $CoV$ with {\normalfont $\text{Re}_0=100$}.}
\end{figure} 
Figure \ref{fig:SE_varyingCoV} shows the performances of the full-rank and the low-rank methods for varying $CoV = \{1\%, 5\%, 10\%\}$. We use Algorithm \ref{alg:adaptive_proc} with 4 Picard steps, followed by several Newton steps until convergence. For $CoV=\{1\%, 5\%\}$, one Newton step is required for the convergence and, for $CoV=10\%$, two Newton steps are required. %We use stringent truncation tolerances in lrGMRES for solving the problems with larger $CoV$ (i.e., we set $\rho_{\text{trunc}} = 10^{-1.5}$ for $CoV = \{1\%, 5\%, 10\%\}$, respectively).
Figure \ref{fig:cost_SE_varyingCoV} shows the computational costs. For $CoV=\{1\%, 5\%\}$, the computational benefits of using the low-rank approximation methods are pronounced whereas, for $CoV=10\%$, the performances of the two approaches are essentially the same for shorter correlation lengths. Indeed, for higher $CoV$, the ranks of solutions $\bar u$ (see Figures \ref{fig:rank_SE_varyingCoV_ux}--\ref{fig:rank_SE_varyingCoV_p}) as well as updates $\delta \bar u^k$ at Newton steps become close to the full rank ($n_\xi = 56$).%, and lrGMRES(20) requires more cycles to generate an approximate satisfying $\truncgmres{}^k$ than GMRES(20) does.  

\begin{figure}[!tb]
	\label{fig:SE_varyingRe}
	\centering
	\subfloat[Computational cost of full-rank computation and low-rank approximation]  {
	\includegraphics[angle=0, scale=0.60]{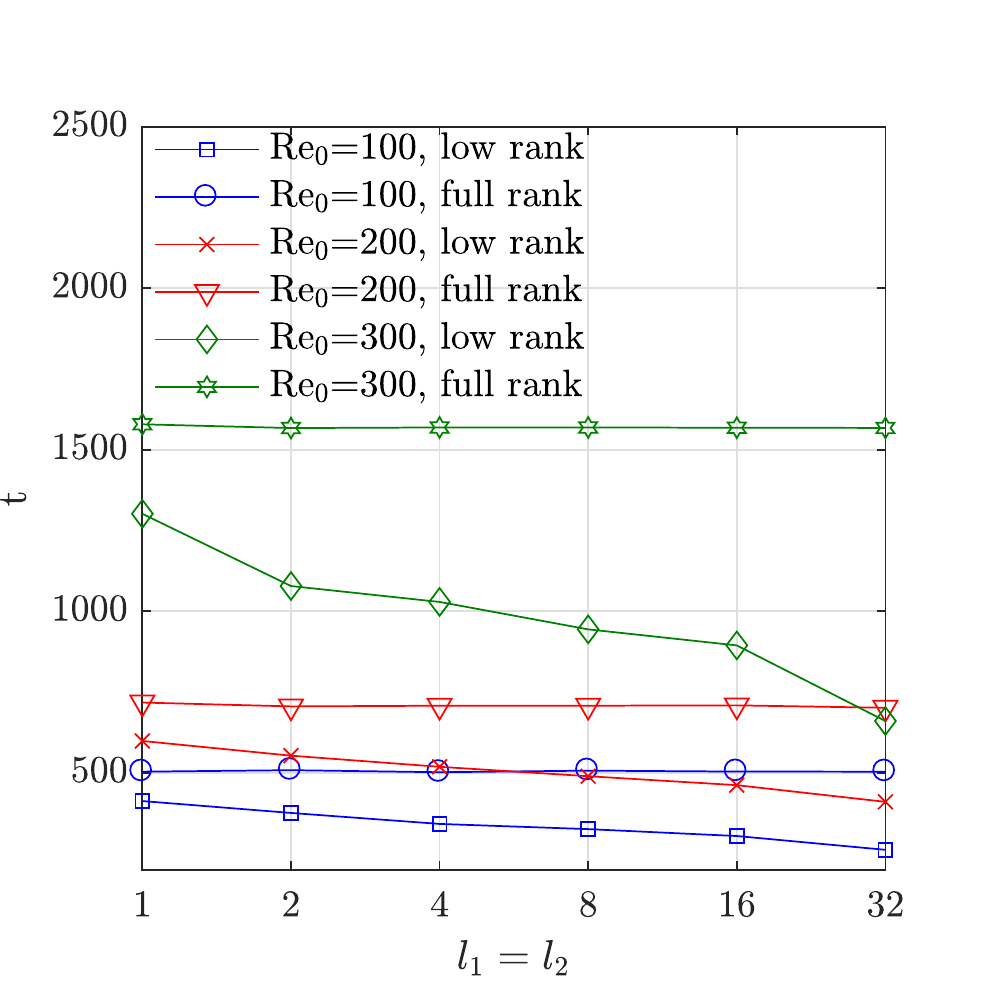}
	\label{fig:cost_SE_varyingRe}
	}\hspace{3mm}	
	\subfloat[Ranks of the low-rank approximate solutions $u_x$]  {
	\includegraphics[angle=0, scale=0.60]{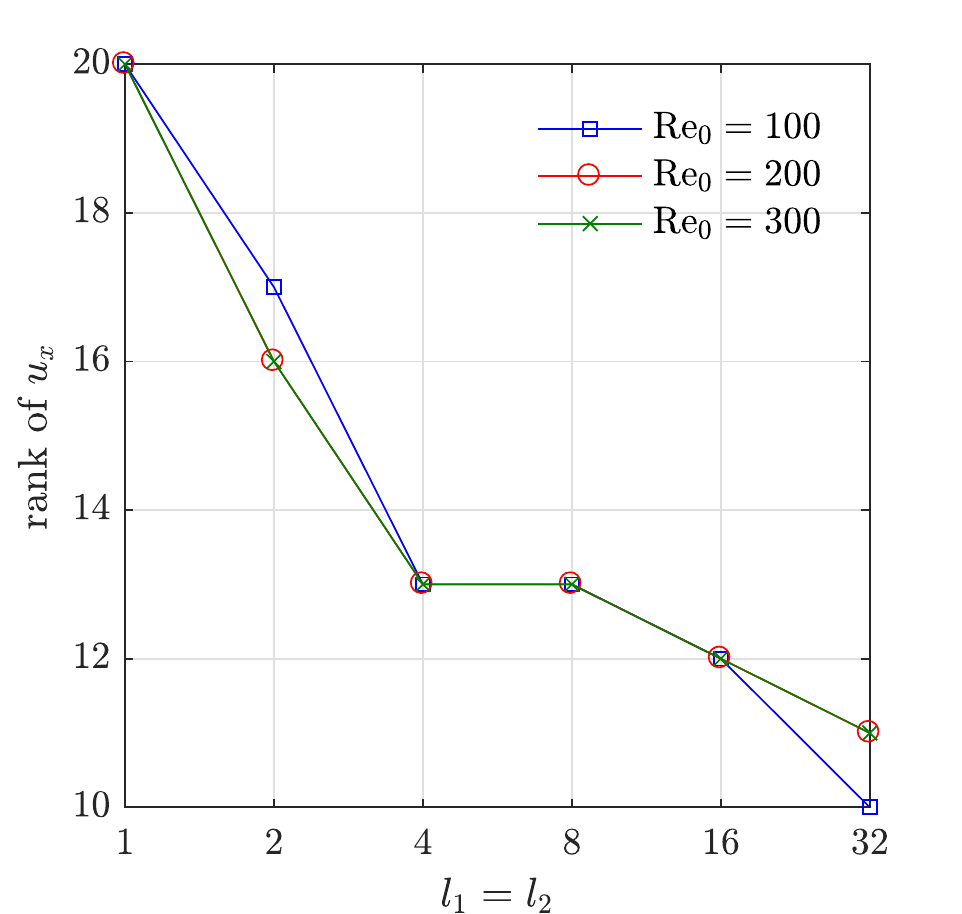}
	\label{fig:rank_SE_varyingRe_ux}
	}\\
	\subfloat[Ranks of the low-rank approximate solutions $u_y$]  {
	\includegraphics[angle=0, scale=0.60]{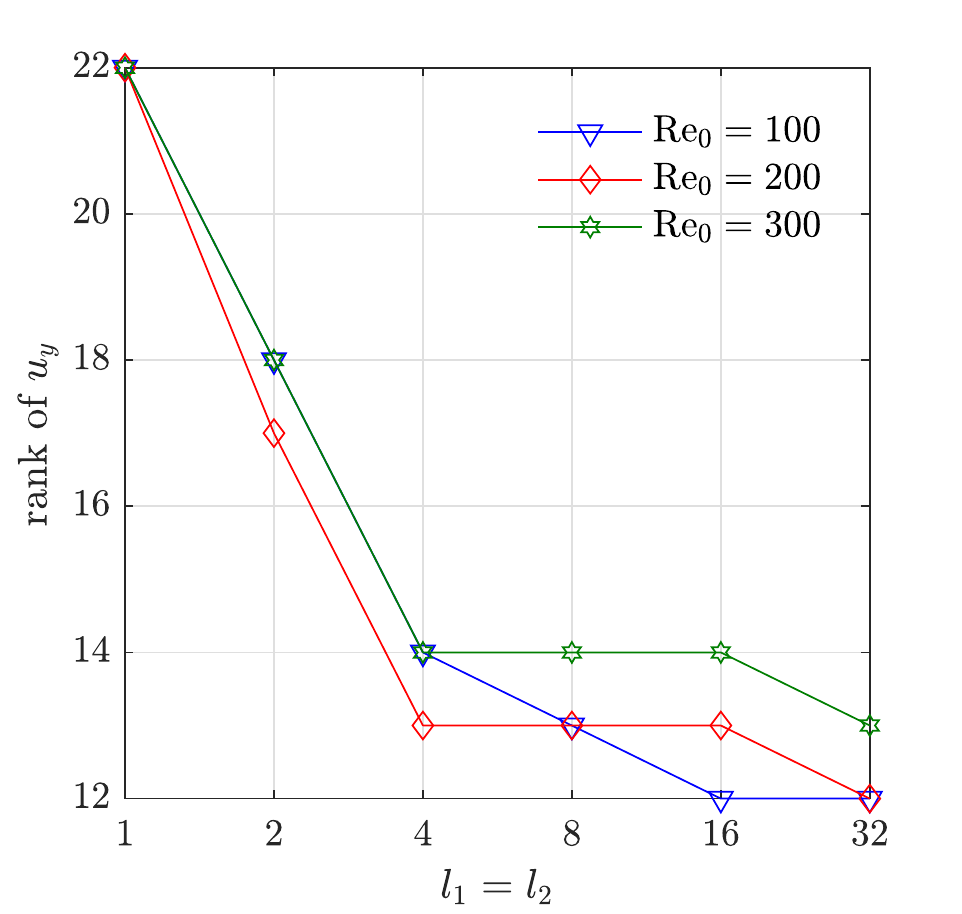}
	\label{fig:rank_SE_varyingRe_uy}
	}\hspace{3mm}	
	\subfloat[Ranks of the low-rank approximate solutions $p$]  {
	\includegraphics[angle=0, scale=0.60]{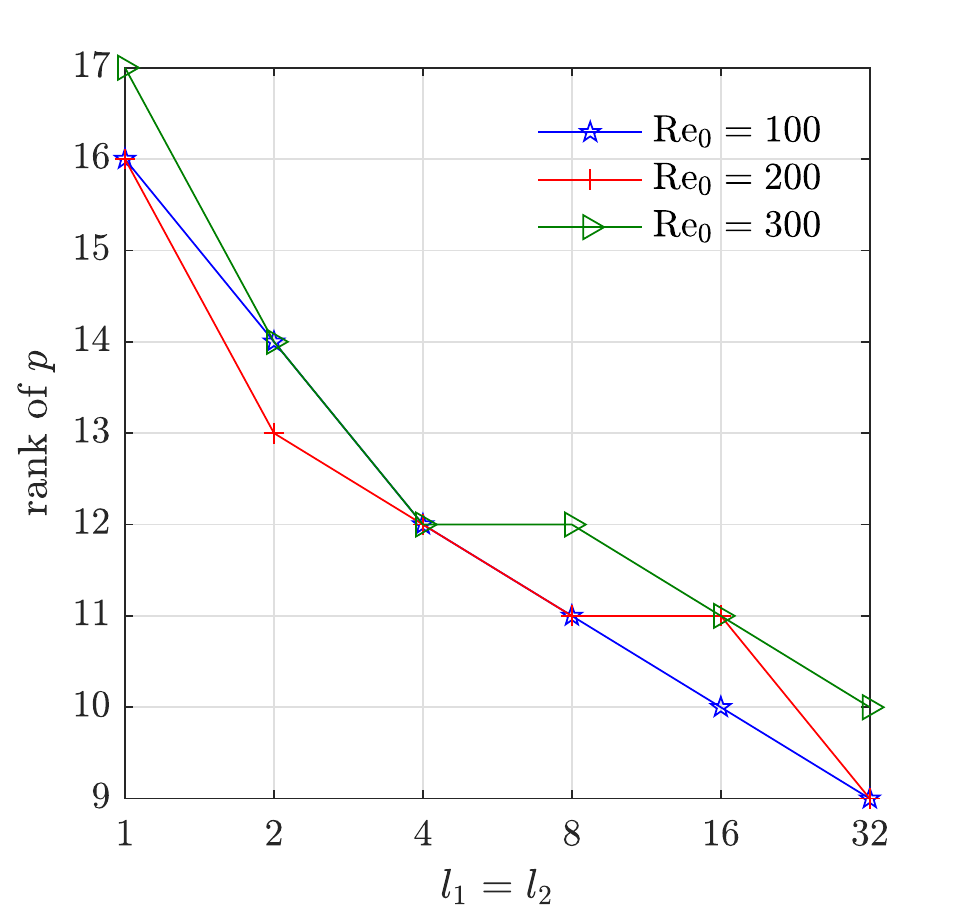}
	\label{fig:rank_SE_varyingRe_p}
	}
	\caption{Computational costs and ranks for varying correlation lengths and varying {\normalfont$\text{Re}_0$}}
\end{figure}

Lastly, we study the benchmark problems with varying mean viscosity with SE covariance kernel and $CoV=1\%$. As the mean viscosity decreases, $\text{Re}_0$ grows, and the nonlinear problem tends to become harder to solve, and for the larger Reynolds numbers $\text{Re}_0 = 200$ or $300$, we use more Picard steps (5 or 6, respectively) before switching to Newton's method.% and use stringent truncation tolerance in lrGMRES by setting $\rho_{\text{trunc}} = 10^{-1.5}$.

Figure \ref{fig:SE_varyingRe} shows the performances of the low-rank methods for varying Reynolds number, $\text{Re}_0=\{100,200,300\}$. For $\text{Re}_0=200$, after 5 Picard steps, one Newton step leads to convergence (and 6 Picard steps and one Newton step for $\text{Re}_0=300$). As the figures \ref{fig:rank_SE_varyingRe_ux}--\ref{fig:rank_SE_varyingRe_p} show, the ranks of the solutions increase slightly as the Reynolds number becomes larger and, thus, for all $\text{Re}_0$ tested here, the low-rank method demonstrates notable computational savings (with $CoV = 1\%$). Note that overall computational costs in Figure \ref{fig:cost_SE_varyingRe} increase as the Reynolds number becomes larger because (1) the number of nonlinear steps required to converge increases as the Reynolds number increases and (2) to solve each linearized systems, typically more lrGMRES cycles are required for the problems with higher Reynolds number.

%% file: conclusion.tex
\section{Conclusion} \label{sec:conclusion}
In this study, we have developed the inexact low-rank nonlinear iteration for the solutions of the Navier--Stoke equations with uncertain viscosity in the stochastic Galerkin context. At each step of the nonlinear iteration, the solution of the linear system is inexpensively approximated in low rank using the tensor variant of the GMRES method. We examined the effect of the truncation to an accuracy of the low-rank approximate solutions by comparing those solutions to the ones computed using exact, inexact nonlinear iterations in full rank and the Monte Carlo method. Then we explored the efficiency of the proposed method with a set of benchmark problems for various settings of uncertain viscosity. The numerical experiments demonstrated that the low-rank nonlinear iteration achieved significant computational savings for the problems with smaller $CoV$ and larger correlation lengths. The experiments also showed that the mean Reynolds number does not significantly affect the rank of the solution and the low-rank nonlinear iteration achieves computational savings for varying Reynolds number for small $CoV$ and large correlation lengths.